\input amstex.tex

\input amsppt.sty

\TagsAsMath

\magnification=1200

\hsize=5.0in\vsize=7.0in

\hoffset=0.2in\voffset=0cm

\nonstopmode

\document

\def\la{\langle}
\def\ra{\rangle}

\input amstex.tex
\input amsppt.sty
\TagsAsMath \NoRunningHeads \magnification=1200
\hsize=5.0in\vsize=7.0in \hoffset=0.2in\voffset=0cm \nonstopmode

\document

\topmatter

\title{On asymptotic stability   of  standing waves of discrete
 Schr\"odinger equation  in $\Bbb Z$ }
\endtitle

\author
Scipio Cuccagna and  Mirko Tarulli
\endauthor

\address
DISMI University of Modena and Reggio Emilia, via Amendola 2,
Padiglione Morselli, Reggio Emilia 42100 Italy\endaddress \email
cuccagna.scipio\@unimore.it \endemail

\address
Dipartimento di Matematica ``L. Tonelli'', University of Pisa, Largo
Pontecorvo 5, 56127 Pisa, Italy\endaddress \email
tarulli\@mail.dm.unipi.it \endemail

\abstract We prove an analogue of a classical asymptotic stability
result of standing waves of the  Schr\"odinger equation originating
in work by Soffer and Weinstein. Specifically, our result is a
transposition on the lattice $\Bbb Z$ of a result by Mizumachi
\cite{M1} and it involves a discrete Schr\"odinger operator
$H=-\Delta +q$.
 The decay rates on the potential are less stringent than
in \cite{M1},  since we require $q\in \ell ^{1,1 }$. We also prove
 $|e^{itH}(n,m)|\le C \langle t \rangle ^{-1/3}$   for a fixed  $C$
  requiring,
in analogy to Goldberg \& Schlag \cite{GSc},  only $q\in \ell
^{1,1}$ if $H$ has no resonances and $q\in \ell ^{1,2}$ if it has
resonances. In this way  we ease the    hypotheses on $H$ contained
in Pelinovsky \& Stefanov \cite{PS}, which have a similar dispersion
estimate.

\endabstract

\endtopmatter

 \head \S 1 Introduction \endhead
 We
consider the   discrete Laplacian $\Delta$ in $\Bbb Z$ defined by
$$(\Delta u)(n)= u(n+1)+u(n-1)-2u(n) .$$
 In $\ell ^2 (\Bbb
Z)$ we have for the spectrum $\sigma (-\Delta
 )=[0,4 ]$. Let for $\langle n \rangle=\sqrt{1+n^2}$
$$\aligned & \ell ^{p,\sigma}(\Bbb Z)=\{ u=\{ u_n\} :
\| u\|  _{\ell ^{p,\sigma}} ^p=\sum _{n\in \Bbb Z} \langle n \rangle
^{p\sigma }|u(n)|^p<\infty \} \text{ for $p\in [1,\infty )$}\\& \ell
^{\infty ,\sigma}(\Bbb Z)=\{ u=\{ u(n)\} : \| u\|  _{\ell
^{\infty,\sigma}}
 =\sup _{n\in \Bbb Z} \langle n \rangle ^{ \sigma }|u(n)| <\infty
\}  .\endaligned $$ We will denote $\ell ^{p,\sigma}(\Bbb Z)$ with
$\ell ^{p,\sigma}$. We will set $\ell ^{p  } =\ell ^{p,0}.$ We will
write $\ell ^{p,\sigma}(\Bbb Z,\Bbb R)$ when we restrict to
functions such that $u_n\in \Bbb R$ for all $n$.

We consider a potential $q=\{ q(n), n\in \Bbb Z \}$ with $q(n)\in
\Bbb R$ for all $n$.
 We consider the discrete Schr\"odinger operator $H$
$$(Hu)(n)=-( \Delta  u)(n)+q(n)u(n).\tag 1.1$$

We   assume:

{\item{(H1)}} $q\in \ell ^{1,1} $.

{\item{(H2)}} $H$ is generic, in the sense of Lemma 5.3.

{\item{(H3)}} $\sigma _d(H)$ consists of exactly one eigenvalue
$-E_0$, with $-E_0\not \in [0,4].$

\bigskip By Lemma 5.3 below $\dim \ker (H+E_0)\le 1 $ in $\ell ^2$. We denote by $\varphi _0(n)$ a
generator of $\ker (H+E_0) $ normalized so that $\| \varphi _0\|
_{\ell ^2}=1.$   Consider now   the
 discrete nonlinear Schr\"odinger equation (DNLS)

$$i\partial _t u (t,n) -  (Hu )(t,n) +     |u(t,n)|^ 6u(t,n)=0
    .\tag 1.2$$
We  look at a particular family of solutions $e^{ i\omega t}\phi
_\omega$ of (1.2), or equivalently of
$$    (Hu )( n) -     |u( n)|^ 6u( n)=-\omega  u ( n)
.\tag 1.3$$ By standard bifurcation arguments we have the following
result, see Appendix A:

\proclaim{Lemma 1.1} Assume (H1)--(H3). There is a family $\omega
\to \phi _\omega $ of standing waves   solving (1.3) with the
following properties. For any $\sigma \ge 0$ there is an $\eta
>0$ such that $\omega \to \phi _\omega  $ belongs to
  $C^\omega  (] E_0, E_0+\eta [, \ell ^{2,\sigma})
  \cap C^0 ([E_0, E_0+\eta  [, \ell ^{2,\sigma})$.
  We have $\phi _\omega (n)
\in \Bbb R$ for any $n$ and  there are   fixed $a>0$ and $C
>0$ such that $|\phi _\omega (n)|\le C e^{-a|n|}$.
As $\omega  \to E_0$   we have in $C^\infty  (] E_0, E_0+\eta [,
\ell ^{2,\sigma})
  \cap C^0 ([E_0, E_0+\eta  [, \ell ^{2,\sigma})$
   the expansion
$$\phi _\omega =(\omega - E _0)^{\frac 1 {6}}  \|
\varphi _0  \| _{\ell ^8} ^{-\frac {4}{3}}(\varphi _0+ O(\omega - E
_0)) .\tag 1.4$$
\endproclaim

  The main aim of this paper
is  the following asymptotic stability result:

\proclaim{Theorem 1.2} Consider in Lemma 1.1 $\sigma >0$ large and
$\eta >0$ small.  Assume (H1)--(H3). For any $\omega _0\in ] E_0,
E_0+\eta [$  there exist  an $\epsilon_0>0$ and a $C>0$ such that if
we pick $u_0\in \ell ^2$ with $  \|u_0- \phi_ {\omega _0} \|_{\ell
^2}<\epsilon <\epsilon _0,$ then there exist $\omega_+\in ( E_0, E_0
+\eta _0 ) $, $\Theta\in C ^1 ( \Bbb R)$ and $ u_+\in
 \ell ^2 $    with $|\omega _+-\omega _0|+\| u_+\| _{\ell ^2} \le C\epsilon  $ such that
 if $u(t,n)$ is the corresponding solution of (1.2) with
 $u(0,n)=u_0(n)$, then
$$\aligned &
\lim_{t\to\infty}\|u(t )-e^{i\Theta(t)}\phi_{\omega_{+}} -e^{
it\Delta }u_+\|_{\ell ^2}=0 .
 \endaligned $$

\endproclaim
For   more precise statements see Theorem  4.1 and Lemma 4.4 in \S
4. Theorem 1.2 is related to \cite{SW2}. The series  \cite{SW1-2}
inspired a long list of papers on asymptotic stability of both large
and small ground states in the continuous case, see
\cite{PW,BP1-2,SW3,Wd,C1,TY1-3,Ts,C2,BS,P,RSS,SW4,GNT,S,KS,GS1,
C3,M1-2,GS2,CM,C4-5,KZ,CT,CV2}. We refer for a discussion of the
state of the art to the introductions in \cite{CM,C5}.  In the case
of the lattice $\Bbb Z$ the results on dispersion in
\cite{SK,KKK,PS} are enough to develop an exactly analogous theory.
The fact that the dispersion rate of $e^{it\Delta
 } $ in $\Bbb Z$ is $\langle t\rangle ^{-1/3}$, \cite{SK}, instead of
$ t^{-1/2}$ in  $\Bbb R$, is analogous, in fact easier, to the
situation in \cite{C3,CV2} which considers operators with potentials
periodic in space. It is possible to develop in   $\Bbb Z$ a theory
completely analogous  to the one in $\Bbb R$, allowing $H$ to have
any finite number of eigenvalues, if, in analogy to the continuous
case, we assume an important nonlinear hypothesis,   the so called
Fermi Golden Rule (FGR) in \cite{CM,C5}. However, the fact that the
spectrum $\sigma
 (-\Delta )$ is $[0,4]$ in $\Bbb Z$ instead of $[0,\infty )$
  in $\Bbb R$, implies that there will be cases in $\Bbb Z$ when
certainly the FGR does not hold. In that case Theorem 1.2 does not
hold any more, see \cite{C6}.

The theory can be developed also for large solitons, following the
framework in \cite{KS}. The fact of the slow $\langle t\rangle
^{-1/3}$ dispersion rate is reflected in the fact that, since we use
Strichartz like estimates, we can prove Theorem 1.2 with
nonlinearity $|u| ^{p-1}u$ for $p\ge 7$, rather than for $p\ge 5$ as
in \cite{M1,C4}. In (1.2) we choose nonlinearity $|u| ^{6}u$ for
definiteness, because what matters is to have a nonlinearity $\beta
(|u|^2)u$ with $\beta (t) \in C^2$ with $\beta (0)=d\beta
(0)/dt=d^2\beta (0)/dt^2=0.$ Notice that the $C^\omega $ regularity
proved in Lemma 1.1, which is used in \cite{C6}, is not necessary
here: all we need is Lemma 1.1 with  $C^\omega $ and $C^\infty$
replaced by $C^2$.

For the proof of the nonlinear estimates and the use of Kato type
smoothing estimates, see Lemmas 3.3-5, we follow  in spirit
\cite{M1} but we also make several simplifications. We recall that
\cite{M1} extends to dimension 1 a strengthening of the result by
\cite{SW2} due to \cite{GNT}.One of the main ingredients in
\cite{GNT} is the endpoint Strichartz estimate, involving spaces
$L^2_tL^{\frac{2D}{D-2}}$,  for $D\ge 3$. Here the key point is the
$L^2 $ norm in time. The main point in \cite{M1} is the search of a
surrogate of the endpoint Strichartz estimate for $D=1$. In
\cite{M1} there are various new interesting smoothing estimates
involving the $L^2 $ norm in time. Here  we point out that,   to
show
  asymptotic stability, the classical smoothing estimates
in \cite{K} are sufficient.
 As a consequence,   following our argument, it is
 possible to prove Theorem 2 \cite{M1}, the stability result
 in \cite{M1}, assuming  the
 decay condition $V\in L^{1,1}(\Bbb R)$ for the potential $V(x)$.
Kato \cite{K}  smoothing  in this paper follows from simple
estimates on the resolvent $R_H(z)$. Lemma 3.5 below requires an
extension to Birman-Solomjak spaces, proved in \cite{CV2}, of a
result by Christ \& Kieselev \cite{CK}, see Lemma 3.1 \cite{SmS}.

A substantial part of this paper is dedicated in reproving the main
result in \cite{PS}. We first develop some theory of Jost functions
for $H$,  along lines very similar to the first part of \cite{DT}.
We then follow   very closely the treatment of dispersion in low
energies for the continuous case in 1D contained in \cite{GSc}. In
this way we prove:

\proclaim{Theorem 1.3} Assume $q\in \ell ^{1,1}$ if $H$ has no
resonances in $0$ and $4$, i.e. $H$ satisfies (H2), and $q\in \ell
^{1,2}$  if 0 or 4 is a resonance, i.e. $H$ does not satisfy (H2).
Then for $P_c(H)$ the projection on the continuous spectrum, we have

$$\| P_c(H)e^{itH} : \ell ^1(\Bbb Z) \to \ell ^\infty (\Bbb Z)\| \le
C \langle t \rangle ^{-1/3}  \text{ for a fixed $C>0$.}$$
\endproclaim
  Theorem  1.3 (which is restated as
Theorem 5.10 and proved in \S 5) is similar to the main result in
\cite{PS}. However here we allow resonances and improve the decay
rates in \cite{PS}, which state their result for $q\in \ell
^{1,\sigma }(\Bbb Z)$ with $\sigma
>4$ and only for the non resonant case.

Few weeks after we produced this work, Kevrekidis {\it et al.}
\cite{KPS} produced independently a result similar to our Theorem
1.2 and Theorem 4.1. They assume $ q\in \ell ^{1,\sigma }$ for
$\sigma
>5$  instead of $\sigma =1$.
 They prove   weaker Strichartz estimates than the ones in Theorem 4.1
  since they do not exploit
 the fact that the upper bound in the dispersive estimate
 in Theorem 1.3 is not
 singular for $t=0$. They do not discuss asymptotic
 completeness. They reproduce the  smoothing argument
 in \cite{M1}. In doing so they avoid  some of the difficulties  in  \cite{M1}
 thanks to the fact that
 energy and $\ell ^2$ norms coincide in the discrete setting.
Nonetheless, this treatment of   smoothing, as well as reliance on
\cite{PS}, is responsible for the more restrictive hypotheses on the
decay of $q(n)$ required in \cite{KPS}.

\bigskip
We end with some notation. Given an operator $A$ we set
$R_A(z)=(A-z)^{-1}$ its resolvent. Given a function $f(\theta )$ for
$\theta \in [-\pi , \pi ]$, we denote by $\widehat{f}(n)$ or by $f
^{\wedge}(n)$ its $n$-th Fourier coefficient, writing Fourier series
in terms of exponentials. We set $f ^{\vee}(n)=f ^{\wedge}(-n)$. We
set $\dot f =\partial _\theta f $; $n^\pm=\max(\pm n, 0 )$. We set
$$\eta (n)=\sum _{m=n}^{\infty} |q(m)|  \text{ and }
\gamma (n)=  \sum _{m=n}^\infty (m-n) |q(m)| .$$ We will denote by
$\Cal S(\Bbb Z)$  the set of functions $f(n)$
 rapidly decreasing as $|n|\nearrow \infty$. We will denote by
 $\Cal S(\Bbb R\times \Bbb Z)$  the set of functions $f(t, n)$
 rapidly decreasing as $(t,n) $ diverges along with all the
 derivatives $\partial _t^af(t, n)$ for $a\in \Bbb N$.
 Given two Banach spaces $X$ and $Y$, $B(X,Y)$  will be the space of
 bounded linear operators defined in $X$ with values in $Y$. By
 $\Bbb Z _{\ge a}$  we mean the subset of $\Bbb Z$ formed by
 elements $\ge a$. For $x\in \Bbb R$ the integer part $[x]\in \Bbb
 Z$ is
 defined by $ [x]\le x <
 [x]+1$.

\head \S 2 Linearization, modulation and set up \endhead

By standard arguments it is possible to prove:

 \proclaim{Lemma 2.1(Global well posedness)} The DNLS (1.2)
 is globally well posed, in the sense that any initial value problem
 $u(0,n)=u_0(n)$ with $u_0\in \ell ^2$ admits exactly one solution
 $u(t)\in C^\infty (\Bbb R, \ell ^2)$. The correspondence $u_0\to u(t)$
 defines a continuous map $\ell ^2\to C^\infty ([T_1,T_2], \ell ^2)$
 for any bounded interval $[T_1,T_2] $.

\endproclaim
By an elementary and standard  implicit function theorem argument,
which we skip,  it is possible to prove the following standard
lemma:

\proclaim{Lemma 2.2 (Coordinates near standing waves)} Fix $\omega
_0$ close to $E_0$. Then there are  an $\epsilon _0>0$ and a $C_0>0$
such that any $u_0\in \ell ^2$ with $\| u_0-\phi _{\omega _0}\|
_{\ell ^2} \le \epsilon \le \epsilon _0$ can be written in a unique
way in the form $  u_0=e^{i\gamma (0)} (\phi _{\omega (0)} +r(0))$
with  $|\omega _0-\omega (0)|+|\gamma (0)|+\| r(0)\| _{\ell ^2} \le
C_0\epsilon $ and with $\langle \Re r(0), \phi _{\omega (0)}\rangle
=\langle \Im r(0), \partial _\omega \phi _{\omega (0)}\rangle =0.$
The correspondence $u_0\to (\gamma (0),\omega (0),r(0))$ is a smooth
diffeomorphism.

\endproclaim
Consider the initial datum $u_0(n)$ which we will suppose close to
$\phi _{\omega _0}$. For some $T>0$ and for $0\le t\le T$ the
corresponding solution $u(t,n) $ can be written as $$ u(t,n) = e^{i
\Theta(t)} (\phi _{\omega (t)} (n)+ r(t,n))  \text{ where
$\Theta(t)=\int _0^t\omega (s) ds + \gamma (t)$}.\tag 2.1$$ with
$(\gamma (t),\omega (t), r(t ))$ with $C^\infty $ dependence in $t$
and such that $\langle \Re r(t), \phi _{\omega (t)}\rangle =\langle
\Im r(t), \partial _\omega \phi _{\omega (t)}\rangle =0 $  for $0\le
t\le T$. When we plough the ansatz in (1.2) we obtain

$$\aligned &
  i \partial _t r (t,n)  =
 (H r ) (t,n)  +\omega (t) r(t,n) -4
  \phi _{\omega (t)} ^6(n) r(t,n)   -3
 \phi _{\omega (t)} ^6(n)  \overline{  r }(t,n) \\& + \dot \gamma (t) \phi
_{\omega (t)}(n) -i\dot \omega (t)
\partial _\omega \phi   _{\omega (t)}(n)
+ \dot \gamma (t) r(t,n)
 +
  N (r(t,n)     ) \endaligned \tag 2.2
$$
for $N (r (t,n)    )=O(r^2(t,n))$. In particular
$$| N (r (t,n)    )|\le c_0 (| \phi _\omega ^5(n) r^2(n)|+ |   r
(n)|^7)\tag 2.3$$ for a fixed $c_0$. The condition $\langle \Re
r(t), \phi _{\omega (t)}\rangle =\langle \Im r(t), \partial _\omega
\phi _{\omega (t)}\rangle =0.$ yields
$$ \Cal A(t)  \left [  \matrix
\dot \omega  \\  \dot \gamma
\endmatrix \right ]  = \left [  \matrix
\langle \Im N ,\phi _{\omega (t)}\rangle \\  \langle \Re N ,
\partial _\omega \phi _{\omega (t)}\rangle
\endmatrix \right ] \tag 2.4$$

$$\text{with } \Cal A(t) = \left [  \matrix
 \frac{1}{2}\partial _\omega \| \phi _\omega \| _{\ell ^2} ^2
-\langle \Re r ,\partial _\omega \phi _{\omega  } \rangle & \langle
\Im r , \phi _{\omega  } \rangle
 \\ \langle \Im  r ,\partial _\omega ^2 \phi _{\omega  } \rangle  &
 -\frac{1}{2}\partial _\omega \| \phi _\omega \| _{\ell ^2} ^2
 -\langle \Re r ,\partial _\omega \phi _{\omega  } \rangle
\endmatrix \right ] .\tag 2.5$$
In order to prove Theorem 1.2 we follow a scheme introduced by
Mizumachi \cite{M1}, in particular we will need a number of
dispersive estimates on $e^{-itH}$. Specifically in Theorem 5.10 we
prove:

 \head \S 3 Spacetime estimates for $H $ \endhead
We list a number of linear estimates needed in the stability
argument.  Given an operator $H$   as in (1.1)
 we will denote by $P_d(H)$ the spectral projection on the discrete
spectrum of $H$ and we will set $P_c(H)=1- P_d(H)$. We set $\ell
^2_c(H)=P_c(H) \ell ^2$.
 The first
result, due to Pelinovsky \& Stefanov \cite {PS} but which we
strengthen,  is Theorem 1.3, see Theorem 5.10 below. Our next step
are the   Strichartz estimates. Here we follow an idea in
\cite{CV1}, for much earlier discussion see   \cite{GV}. For every
$1\leq p, q\leq \infty$ we introduce the Birman-Solomjak spaces
$$
\ell ^p({   \Bbb Z}, L^q_t[n,n+1])\equiv \left \{f\in
L^q_{loc}({\Bbb R}) \hbox{ s.t. } \{\|f\|_{L^q[n, n+1]} \}_{n\in
{\Bbb Z}} \in \ell ^p({\Bbb Z})\right \}, $$  endowed with the
  norms
$$\aligned &    \| f \|_{\ell ^p({  \Bbb Z}, L^q_t[n,n+1])} ^p\equiv
\sum_{n\in \Bbb Z}   \| f \| _{L^q_t[n,n+1] }^{p}   \quad  \forall
\quad  1\le p<\infty \text{ and } 1\leq q \leq \infty  \\ & \| f \|
_{\ell ^\infty ({ \Bbb Z}, L^q_t[n,n+1])} \equiv \sup _{n\in \Bbb Z}
\| f \| _{L^q[n,n+1] }.
\endaligned $$
We will say that a pair of numbers $(r,p) $ is admissible if
 $$ 2/{r }+ 1/{p }=1/2 \text{ and }
 (r,p) \in [4, \infty]\times [2, \infty]  .\tag 3.1$$
Then by the standard $TT^\ast$ argument  it is possible to prove
from Theorem 1.3 the following result, whose proof we skip, but see
Lemma 3.2 in \cite{CV2} and the proof in \S 9 \cite{CV2}:

\proclaim{Lemma 3.1 (Strichartz estimates)} Under the hypotheses and
conclusions of Theorem 1.3 there exists a constant $C =C_H $  such
that for every admissible pair $(r,p) $   we have:
$$
  \|e^{itH }P_c(H  ) f \|
    _{\ell ^{\frac 32 r}( \Bbb {Z},L^{\infty}_{t}([n,n+1],
    \ell ^{ p}(\Bbb {Z})))} \leq
C  \|f\|_{\ell ^2(\Bbb {Z})}.
$$
Moreover,  for any two admissible pairs $(r_1,p_1), (r_2,p_2) $   we
have the   estimate
$$ \aligned & \left \| \int _{0}^{t}e^{i(t-s)H }P_c(H  ) g(s )ds\right \| _{\ell ^{\frac{3}{2}r_1}( \Bbb
{Z},L^{\infty }_{t}([n,n+1],\ell ^{ p_1}({\Bbb Z}))} \le \\& \le C
\| g \| _{\ell ^{\left ( \frac{3}{2}r_2\right ) '}( \Bbb
{Z},L^{1}_{t}([n,n+1],\ell  ^{ p_2'}({\Bbb Z}))}.\endaligned
$$\endproclaim
In \S 5 we prove the following Kato smoothness result:

\proclaim{Lemma 3.2} Assume that  $H$ is generic with $q\in \ell
^{1,1}$. For $\tau >1$ there exists $C=C(\tau  )$ such that  for all
$z\in \Bbb C \backslash [0,4]$
$$\| R_H(z )P_c(H  )   \| _{B (\ell ^{2, \tau },
\ell ^{2,-\tau }) }\le  C  .$$ The  following limits are well
defined for any $\lambda \in [0,4]$ in $C^0([0,4],B(\ell ^{2,\tau }
, \ell ^{2,-\tau }  ))$
$$\lim _{\epsilon \to 0 ^+} R_{H }(\lambda \pm i\epsilon )
=R^{\pm}_{H }(\lambda  )   .$$ For any $u\in \ell ^{2,\tau}\cap \ell
^2_c(H)$ we have
$$P_c(H)u=\frac{1}{ 2\pi i}\int_\Bbb R
 (R_{H  }^{+}(\lambda  )-R_{H
}^{-}(\lambda  ))   u d\lambda  =\frac{1}{ 2\pi i}\int_0^4
 (R_{H }^{+}(\lambda  )-R_{H
}^{-}(\lambda  ))   u d\lambda   .$$
\endproclaim

   The  next few lemmas are simplifications of corresponding
   lemmas in \cite{M1}. First of all we skip the smoothing
   estimates in Lemma 4 \cite{M1} and we instead consider the
   following classical consequence of Lemma 3.2, see \cite{K}:

\proclaim{Lemma 3.3} Assume that  $H$ is generic with $q\in \ell
^{1,1}$. Then for $\tau >1$ we have:
   {\item {(a)}}
  for any $f\in \Cal S(\Bbb Z)$ and for $C(\tau )$ the constant of
  Lemma 3.2
$$\align & \| e^{-itH }P_c(H  )f\|
_{\ell  ^{2 , -\tau } L^2_t} \le
 2 \sqrt{ \pi C(\tau )}\|f\|_{\ell ^2};
\endalign $$
  {\item {(b)}}
  for any $g(t,n)\in
 \Cal S (\Bbb R \times \Bbb Z)$
$$ \left\|\int_\Bbb R e^{itH }P_c(H  )g(t,\cdot)dt\right\|_{\ell ^2 } \le
2\sqrt{ \pi C(\tau )}\| g\|_{\ell  ^{2, \tau}L_t^2}.
$$\endproclaim
{\it Proof.} (a) implies (b) by duality. So we focus on (a), which
is a consequence of  \S 5 \cite{K}. Get $g(t,\nu ) \in S(\Bbb
R\times \Bbb Z) $ with $g(t)=P_c(H )g(t)$. By the limiting
absorption principle in Lemma 3.2
$$\aligned & \la  e^{-itH  } f,
g\ra _{L^2_t\ell ^2 }=    \frac{1}{ 2\pi i}\int_0^4
   \left\la
 (R_{H  }^{+}(\lambda  )-R_{H
}^{-}(\lambda  ))   f,  \overline{\widehat{g}}(\lambda
)\right\ra_{\ell ^2} d\lambda .
\endaligned $$
Then from Fubini and  Plancherel   we have
$$\aligned & \big |\la e^{-itH  } f, g\ra _{L^2_t\ell ^2 }\big |
 \le (2\pi)^{-\frac{1}{2}}
\| (R_{H  }^{+}(\lambda  )-R_{H  }^{-}(\lambda )) f\|_{\ell ^{2,
-\tau} L^2_\lambda (0,4)}   \|  {g} \|_{\ell ^{2, \tau}  L^2_t }
 .
  \endaligned \tag 3.2$$
We have
$$\| (R_{H  }^{+}(\lambda  )-R_{H  }^{-}(\lambda )) f\|_{\ell ^{2,
-\tau} L^2_\lambda (0,4)}=\lim _{\varepsilon \searrow 0} \| (R_{H }
(\lambda +i\varepsilon )-R_{H  } (\lambda -i\varepsilon )) f\|_{\ell
^{2, -\tau} L^2_\lambda (0,4)}.
$$
Then we repeat an argument in Lemma 5.5 \cite{K}. For $\Im \mu >0$
set $K(\mu )$ the positive square root of $(2\pi i)^{-1} [R_{H  }
(\mu )- R_{H } (\overline{\mu} )]= \pi ^{-1}(\Im \mu )R_{H  }
(\overline{\mu} ) R_{H  } ( {\mu} )$. Then
$$\aligned &  \int _0^4 \| (R_{H }
(\lambda +i\varepsilon )-R_{H  } (\lambda -i\varepsilon )) f\|_{\ell
^{2, -\tau} } ^2d\lambda =\\& 4\pi ^2 \int _0^4 \| \langle n \rangle
^{-\tau }K(\lambda +i\varepsilon )K (\lambda +i\varepsilon )
f\|_{\ell ^{2 } } ^2d\lambda \le 8\pi ^2 C(\tau )\int _0^4   \|  K
(\lambda +i\varepsilon ) f\|_{\ell ^{2 } } ^2d\lambda \\& \le 8\pi
^2 C(\tau )\| f\|_{\ell ^{2 } } ^2
  \endaligned \tag 3.3$$
with $C(\tau )$ the constant in Lemma 3.2. By (3.3) we get in (3.2)
$$\aligned & \big |\la e^{-itH  } f, g\ra _{L^2_t\ell ^2 }\big |
 \le 2\sqrt{  \pi  C(\tau )}\| f\|_{\ell ^{2 } }  \|  {g} \|_{\ell ^{2, \tau}  L^2_t }
 .
  \endaligned \tag 3.4$$
(3.4) yields claim (a).

 \proclaim{Lemma 3.4} Assume    $H$   generic with $q\in \ell
^{1,1}$. Then for any  $\tau >1$ $\exists$ $C_\tau   $   s.t.
$$\align &  \left\|  \int_0^t e^{-i(t-s)H }
P_c(H)g(s,\cdot)ds\right\|_{\ell ^{2 ,-\tau }L_t^2} \le C_\tau   \|
g\|_{\ell ^{2, \tau }L_t^2}.  \endalign
$$
\endproclaim
{\it Proof.} By Plancherel and H\"older inequalities and by Lemma
3.2

 we have
$$ \aligned &
\| \int _{0}^t e^{-i(t-s)H }P_c( H)g(s,\cdot)ds\|_{ L_{t}^2\ell
^{2,-\tau }} \le \\& \le  \| R_{H  }^+(\lambda )P_c( H)
  \widehat{ \chi }_{[0,+\infty )}\ast _\lambda
   \widehat{ g}(\lambda,x)\|_{L_{\lambda  }^2 \ell ^{2,-\tau } }   \le \\& \le
   \left\| \,
\|   R_{H  }^+ (\lambda )P_c( H) \| _{B( \ell ^{2, \tau }, \ell
^{2,-\tau })} \|
   \widehat{ \chi }_{[0,+\infty )}
   \ast _{\lambda } \widehat{g} (\lambda,x) \|_{\ell ^{2, \tau }}\, \right\|_{L^2_\lambda}
\\ \le & \sup _{\lambda \in \Bbb R }
  \|   R_{H  }^+ (\lambda )P_c( H) \| _{B( \ell ^{2,
\tau }, \ell ^{2,-\tau })}  \| g\|_{L_{t }^2\ell ^{2, \tau }
}\lesssim  \| g\|_{L_{t }^2\ell ^{2, \tau } }.
\endaligned $$

\proclaim{Lemma 3.5}  Assume   $H$   generic with $q\in \ell
^{1,1}$. Then for every $\tau >1$  $\exists$ $C_\tau $ s. t.
$$
\left\|\int_0^t e^{-i(t-s)H }P_c(H  )g(s,\cdot)ds \right\|_{
L_t^\infty \ell  ^2\cap \ell ^{6}( \Bbb {Z},L^{\infty
}_{t}([n,n+1],\ell ^\infty ))} \le C \|g\|_{L_t^2\ell ^{2,\tau }}.$$
\endproclaim
{\it Proof.} For $g(t,\nu )\in S(\Bbb R\times \Bbb Z)$ set
$$T g(t)=\int _0^{+\infty}  e^{-i(t-s)H }P_c(H  )g(s) ds  .$$
Lemma 3.3 (b) implies $f:=\int _0^{+\infty}  e^{isH  }P_c( H
)g(s)ds\in\ell ^2$.  For $(r,p) $ admissible we have
$$ \aligned &  \|T g(t)\|_{\ell ^{\frac 32 r}
( \Bbb {Z},L^{\infty}_{t}([n,n+1],\ell ^{ p} ))}\lesssim  \|
f\|_{\ell ^2} \lesssim  \|g\|_{L^2_t\ell ^{2, \tau}  } .\endaligned
$$  Lemma 3.5  follows from this estimate by an extension  of
the Christ Kieselev  Lemma 3.1 \cite{SmS}  to Birman - Solomjak
spaces. The proof of this  extension is in  \cite{CV2}. Consider two
Banach spaces and $X$ and $Y$   and  $K(s,t)$  continuous function
valued in the space $B(X, Y)$. Let
$$T_K f(t)=\int_{-\infty}^\infty K(t,s)f(s) ds\text{ and }
 \tilde T_K f(t) = \int_{-\infty}^t K(t,s) f(s) ds.$$
 Then we have: \proclaim
{Lemma 3.6} Let $1\leq p, q,r\leq \infty$ be such that $1\leq r<\min
( p,q )\leq \infty$. Assume that there exists $C>0$ such that
$$
\|T_K f\|_{\ell ^q(\Bbb Z, L^p_t([n,n+1],Y))}\leq C
\|f\|_{L^r_t(X)}.$$ Then
$$
\|\tilde T_K f\|_{\ell ^q (\Bbb Z, L^p_t([n,n+1],Y))}\leq C'
\|f\|_{L^r_t(X)} $$ where $C'=C'(C, p,q,r)>0$.
\endproclaim
 In the case $p=q$ the previous lemma follows from \cite{CK},
 while the general case is in \cite{CV2}.

\head \S 4  Proof of Theorem 1.2 \endhead

Our first goal here is to prove the following:

\proclaim{Theorem 4.1} Fix $\omega _0\in ] E_0, E_0+\eta [$ with
$\eta >0$ sufficiently small. Then there exist an $\epsilon_0>0$ and
a $C>0$ such that if we pick $\|u_0- \phi_ {\omega _0} \|_{\ell
^2}<\epsilon <\epsilon _0,$ then   the ansatz (2.1) with $\langle
\Re r(t), \phi _{\omega (t)}\rangle =\langle \Im r(t),
\partial _\omega \phi _{\omega (t)}\rangle =0 $ is valid for all
times  and we have the following inequalities:
$$\align  & \| (\dot \omega , \dot \gamma ) \|
_{ L^1(\Bbb R)\cap L^\infty (\Bbb R)} \le C \epsilon, \tag 1\\&
 \| r(t,n) \|    _{\ell ^{\frac{3}{2}r }(
\Bbb {Z}  ,L^{\infty }_{t}([n,n+1],\ell ^{ p} ))}\le C \epsilon
\, \text{ for all admissible $(r,p)$}. \tag 2
\endalign $$

\endproclaim
{\it Remark}. Notice that (1) implies the existence of $\gamma _\pm
$ and $\omega _\pm $ such that $$\lim _{t\to \pm \infty }   ( \omega
,   \gamma )(t)=( \omega _{\pm} ,   \gamma _{\pm} ).$$ By Lemma 2.2
we conclude that (1) implies that $$|( \omega _{0} ,   0 ) -  (
\omega ,   \gamma )(t) |\le C(\omega _0, C)\epsilon \, \text{ for all
$t$}.$$

\bigskip
{\it Proof of Theorem 4.1.} There are two equivalent ways to prove
results like Theorem 4.1. One way it to  prove estimates (1-2) over
bounded intervals $[-T,T]$ with constants $\epsilon _0,C$
independent of $T$ and then let $T\nearrow \infty $. However one can
reach the same result by assuming that the global space-time
estimates hold for some large constant $C_1$, and then by showing
that the estimates hold also for $C_1/2$. By standard arguments,
this sort of a priori estimates method yields Theorem 4.1. Set
$X_{(r,p)}= \ell ^{\frac{3}{2}r }( \Bbb {Z}  ,L^{\infty
}_{t}([n,n+1],\ell ^{ p} )).$ We will then prove:

\proclaim{Lemma  4.2} Fix $\omega _0\in ] E_0, E_0+\eta [$ and $\eta
>0$   small. For  $D$ with
$$1<D<\eta ^{-1} , \tag 1$$
 there exists an $\epsilon _0=\epsilon _0(D) $
with
$$0<\epsilon _0< (\omega _0-E_0)^7,  \tag 2$$
such that, if for an $u_0$ with $\|u_0- \phi_ {\omega _0} \|_{\ell
^2}<\epsilon <\epsilon _0 $   ansatz (2.1) with $\langle \Re r(t),
\phi _{\omega (t)}\rangle =\langle \Im r(t),
\partial _\omega \phi _{\omega (t)}\rangle =0 $  is
valid for all times and  if we have
$$\align  & \| (\dot \omega , \dot \gamma ) \|
_{ L^1(\Bbb R)\cap L^\infty (\Bbb R)} \le D \epsilon, \tag 3\\&
 \| r(t,n) \|    _{X_{(r,p)}\cap L^2_t\ell ^{2,-2}}\le D \epsilon
\text{ for all admissible $(r,p)$}, \tag 4
\endalign $$
    then  for some fixed $C_0$ and for all
admissible $(r,p)$

$$\align  & \| (\dot \omega , \dot \gamma ) \|
_{ L^1(\Bbb R)\cap L^\infty (\Bbb R)} \le  \epsilon, \tag 5\\&
 \| r(t,n) \|    _{X_{(r,p)}\cap L^2_t\ell ^{2,-2}}\le C_0    \epsilon
  .\tag 6
\endalign $$

\endproclaim
It is enough to prove Lemma 4.2 to obtain Theorem 4.1.     Set
$$P_ d (\omega )=2\frac{\langle \Re \, \cdot \, ,\phi _{\omega } \rangle }
{\partial _\omega \| \phi _{\omega } \| _{\ell ^2}^2} \partial
_\omega \phi _{\omega } + 2i\frac{\langle \Im \, \cdot \, ,\partial
_\omega \phi _{\omega } \rangle } {\partial _\omega \| \phi _{\omega
} \| _{\ell ^2}^2}
  \phi _{\omega }\text{ and } P_c(\omega  )= 1-P_
d(\omega ).$$ Then $P_d(\omega )$ and $P_c(\omega )$ are  well
defined operators in $\ell ^{p,\sigma}$ for all $p\ge 1$ and $\sigma
\in \Bbb R$. We have $P_d(\omega (t))r(t)=0$ for all $t$. By Lemma
1.1 we have
$$\aligned & \partial _\omega \phi _{\omega } =\frac{1}{6}
(\omega - E _0)^{-\frac 5 {6}}  \| \varphi _0  \| _{\ell ^8}
^{-\frac {4}{3}}(\varphi _0+ O(\omega - E _0)) \text{ in any $\ell
^{p,\sigma}$,}
\\& \partial _\omega \| \phi _{\omega }
\| _{\ell ^2}^2 =\frac{1}{3}(\omega - E _0)^{-\frac 2 {3}}  \|
\varphi _0  \| _{\ell ^8} ^{-\frac {8}{3}}(1+ O(\omega - E _0))
.\endaligned $$ With   (1.4) these expansions imply
$$\| P_ d (\omega  ) -P_ d (H  )\|
_{B(\ell ^{p,\sigma}, \ell ^{p,\sigma} )}   \le C_{p,\sigma }
|\omega - E_0|  . $$ Then (1-3) imply $C_{p,\sigma } |\omega (t) -
E_0|  \le 1/2$ for all $t$. So from $r(t)=P_ c (\omega (t) )r(t)$ we
conclude $ \| r( t)\| _{\ell ^{p,\sigma }}/ \| P_c(H) r( t)\| _{\ell
^{p,\sigma }} \in [1/2,2]$. Hence $ \| r \| _{X_{(r,p)}}/  \| P_c(H)
r \| _{X_{(r,p)}}$ and $\| r \| _{L^2_t\ell ^{2,-2}}/
 \| P_c(H) r \| _{L^2_t\ell ^{2,-2}} $ are in $[1/2,2]$. We rewrite (2.2) in the form

$$\aligned &
  i \partial _t r (t,n)  =
 (H r ) (t,n)  +\omega (t) r(t,n)+ \dot \gamma (t) r(t,n) +
 \sum _{j=2}^4g_j,\\&
g_2(t,n)= -4
  \phi _{\omega (t)} ^6(n) r(t,n)   -3
 \phi _{\omega (t)} ^6(n)  \overline{  r }(t,n), \\&
 g_3(t,n)= \dot \gamma (t) \phi
_{\omega (t)}(n) -i\dot \omega (t)
\partial _\omega \phi   _{\omega (t)}(n),
\\& g_4(t,n)=
  N (r(t,n)    ) .\endaligned \tag 4.1
$$
Set $w(t)=  e^{ i\Theta (t)}r(t )$ with $\Theta (t) =\int _0^t
\omega (s) ds+\gamma (t)$. Obviously  $\| P_c(H) r \|
_{X_{(r,p)}}=\| P_c(H) w \| _{X_{(r,p)}}$  $\forall$ $X_{(r,p)}$ and
$\| P_c(H) r \| _{L^2_t\ell ^{2,-2}}=\| P_c(H) w \| _{L^2_t\ell
^{2,-2}}$ . We have
$$\aligned & P_c(H)w(t)=e^{-iHt}P_c(H)w(0)+\sum _{j=2}^4 w_j(t),
 \\&
w_j(t)=-i \int _0^t e^{-iH(t-s) } e^{ i\Theta (s)}P_c(H)g_j (s)
ds.\endaligned\tag 4.2$$ The following implies Lemma 4.2:

\proclaim{Lemma  4.3} Assume the hypotheses of Lemma 4.2.
 Then
we have   for some fixed $C_1$, not dependent on $\omega _0$ or any
other parameter,
 and for all admissible
$(r,p)$

$$\align  & \|  \dot \omega \| _{ L^1(\Bbb R)\cap L^\infty (\Bbb R)}\le C_1D^2
(\omega _0-E_0) ^{5/3} \epsilon ^2\tag 7\\& \|\dot \gamma   \| _{
L^1(\Bbb R)\cap L^\infty (\Bbb R)} \le C_1D^2 (\omega _0-E_0) ^{2/3}
\epsilon ^2 \tag 8\\&
 \| w_2(t,n) \|    _{X_{(r,p)}\cap L^2_t\ell ^{2,-2}}\le C_1 D (\omega
_0 -E_0)\epsilon  , \tag 9\\& \| w_3(t,n) \|    _{X_{(r,p)}\cap
L^2_t\ell ^{2,-2}}\le C_1D^2 (\omega _0-E_0) ^{\frac{5}{6}}
 \epsilon ^2  , \tag 10\\& \|
w_4(t,n) \| _{X_{(r,p)}\cap L^2_t\ell ^{2,-2}}\le C_1D^2 (\omega
_0-E_0)^{\frac{5}{6}}\epsilon ^2  .\tag 11
\endalign $$
\endproclaim
Notice that by Lemmas 3.1-2 we have for fixed constants  $$\|
e^{-itH}P_c(H)w(0)\| _{X_{(r,p)}\cap L^2_t\ell ^{2,-2}} \lesssim \|
w(0)\| _{\ell ^2}\lesssim \epsilon .$$  Then Lemma 4.3  and
hypotheses (1)--(2) in Lemma 4.2 imply for some fixed $C_H$
$$\| P_c(H) w(t) \| _{X_{(r,p)}\cap L^2_t\ell ^{2,-2}}\approx \| w(t)
\| _{X_{(r,p)}\cap L^2_t\ell ^{2,-2}}\le C_H \epsilon .$$ This
yields Lemma 4.2.

\noindent We focus now on Lemma 4.3. We observe that (2)--(3) yield
$ \frac{|\omega -E_0|}{|\omega _0 -E_0|}\in [1/2,3/2]$ by
$$|\omega (t)-\omega _0|\le \|\dot \omega \| _{L^1}+|\omega (0)-\omega _0|\le
(D  +C_0)\epsilon <|\omega _0 -E_0|/2,$$ where we used also Lemma
2.2. Now we prove
 (7-11).
We start from (7). By the hypotheses in Lemma  4.2, $\Cal
A(t)\approx
\partial _\omega \|\phi _{\omega } \| ^{2}_{\ell ^2} 2^{-1}\text{diag}(1,-1)+o(\partial
_\omega \|\phi _{\omega } \| ^{2}_{\ell ^2}) $ in $GL(2)$, with
$\Cal A$ the matrix in (2.5).  Then $\| \Cal A^{-1}(t)\| _{GL(2)}
\lesssim  (\omega _0-E_0) ^{2/3}$. Thus by (2.3)--(2.4), by (2) and
by (4)
$$\aligned &  |\dot \omega (t)| \lesssim (\omega _0-E_0) ^{\frac 23}
\| \phi _\omega \| _{ \ell ^{\infty  }} \left (  \| \phi _\omega \|
_{ \ell ^{\infty  }} ^5 \| w (t) \| _{ \ell ^{2 } }^{2} +\| w (t) \|
_{ \ell ^{2 } }^{7}\right )
 ,
 \\& \|\dot \omega \| _{L^1_t }\lesssim (\omega _0-E_0) ^{\frac 23}
\| \phi _\omega \| _{ L^\infty _t\ell ^{\infty , 4 }} \left(  \|
\phi _\omega \| _{ L^\infty _t\ell ^{\infty   }} ^5 \| w \| _{L^2_t
\ell ^{2,-2}  }^{2} + \| w \| _{L^\infty _t \ell ^{2 }  }  \| w \|
_{L^6_t \ell ^{\infty }  }^{6}\right ) .\endaligned \tag 12$$ So $\|
\dot \omega \| _{L^1  \cap L^\infty} \lesssim  D^2 (\omega _0-E_0)
^{\frac 5 3} \epsilon ^2$. Similarly $\| \dot \gamma \| _{L^1  \cap
L^\infty} \lesssim  D^2 (\omega _0-E_0) ^{\frac 2 3} \epsilon ^2$ by

$$\aligned &  |\dot \gamma (t)| \lesssim (\omega _0-E_0) ^{\frac 23}
\| \partial _\omega \phi _\omega \| _{ \ell ^{\infty   }}  \left (
\| \phi _\omega \| _{ \ell ^{\infty  }} ^5 \| w (t) \| _{ \ell ^{2 }
}^{2} +\| w (t) \| _{ \ell ^{2 } }^{7}\right ),
 \\& \|\dot \gamma \| _{L^1_t }\lesssim (\omega _0-E_0) ^{\frac 23}
\| \partial _\omega \phi _\omega \| _{ L^\infty _t\ell ^{\infty , 4
}}
   \left(  \|
\phi _\omega \| _{ L^\infty _t\ell ^{\infty   }} ^5 \| w \| _{L^2_t
\ell ^{2,-2}  }^{2} + \| w \| _{L^\infty _t \ell ^{2 }  }  \| w \|
_{L^6_t \ell ^{\infty }  }^{6}\right ).\endaligned  \tag 13$$ We
have for $j=2,3,4$

$$\| w_j(t) \| _{X_{(r,p)}\cap L^2_t\ell  ^{2,-2} }
\le \left \| \int _0^t e^{-iH(t-s) } e^{ i\Theta (s)}g_j (s) ds
\right \| _{X_{(r,p)}+  L^2_t\ell  ^{2,-2}}.\tag 14
$$
For $j=2,3$ the latter   $\lesssim \| g_j \| _{\ell ^{2,  2}L^2_t}$
by Lemmas 3.4-5. We obtain
$$\aligned & \| g_2 \| _{\ell ^{2,
2}L^2_t} \le 7\| \phi _{\omega } \| _{\ell ^{\infty , 2}L^\infty
_t}^6 \| w \| _{\ell ^{2, - 2 }L^2_t}
\\& \lesssim
(\omega _0-E_0)   \| w \| _{\ell ^{2, - 2 }L^2_t} \le D (\omega
_0-E_0)   \epsilon \le  \epsilon .\endaligned \tag 4.3$$ Moreover we
get $ \| g_3 \| _{\ell ^{2, 2}L^2_t} \le \| \dot \omega \partial
_\omega \phi _{\omega } \| _{\ell ^{2, 2}L^2_t} + \| \dot \gamma
\phi _{\omega } \| _{\ell ^{2, 2}L^2_t}.$ Then by (12)--(13) we have
$$ \aligned & \| g_3 \| _{\ell ^{2,
2}L^2_t} \le \| \dot \omega \| _{L^2_t} \|
\partial _\omega \phi _{\omega } \| _{L^\infty _t\ell
 ^{2, 2} } +
\| \dot \gamma \| _{L^2_t} \|
  \phi _{\omega } \| _{L^\infty _t\ell
 ^{2, 2} }\\& \lesssim D^2 (\omega _0-E_0) ^{\frac{5}{6}}
 \epsilon ^2 <\epsilon .
\endaligned \tag 4.4$$
Finally we need to bound $\| w_4\|  _{X_{(r,p)}\cap L^2_t\ell
^{2,-2} }$. We split $|g_4|\lesssim  |\phi _\omega |^5|w|^2 + |w|^7
$. Correspondingly we have $w_4= w_{4,1}+w_{4,2}$ with, by the above
arguments

$$\| w_{4,1}\|  _{X_{(r,p)}\cap L^2_t\ell ^{2,-2} }\lesssim
(\omega _0-E_0)^{\frac{5}{6}} \| w \|  _{  L^2_t\ell ^{2,-2} }^2\le
D^2 (\omega _0-E_0)^{\frac{5}{6}}\epsilon ^2 .\tag 4.5
$$ By Lemma 3.1, $$ \aligned & \| w_{4,2}\| _{X_{(r,p)} } \le   C_0
\| w^7\| _{L^1_t\ell ^2} \le C_0\| w \| _{ L^{ \infty }_{t} \ell
^{2 } } \| w   \| ^6 _{ L^{6}_{t}\ell ^{\infty }  }
\\&  \le C_0 \| w \| _{ L^{ \infty }_{t} \ell  ^{2  } } \| w\| ^6 _{
\ell ^6 (\Bbb {Z},L^\infty _t[n, n+1]),\ell ^\infty ))} \le C_0D^7
\epsilon ^7.
\endaligned \tag 4.6$$
By (b) Lemma 3.3 and by the argument in (4.6)
$$\aligned & \| w_{4,2}\|  _{  L^2_t\ell  ^{2,-2} }
\le  \int _0^\infty  \|
 e^{-i(t-s) H }
P_c(H)O(w ^7)  (s)\|  _{  L^2_t\ell  ^{2,- 2} } ds    \\& \le C_H
\int _0^\infty \| O(w ^7) (s) \| _{ \ell ^2}ds\lesssim D^7 \epsilon
^7 .
\endaligned \tag 4.7$$
Then (4.5)--(4.7) yield estimate (11) in Lemma 4.3 by hypothesis (2)
in  Lemma 4.2. This completes the proof of Lemma 4.3.

The following standard lemma yields the asymptotic flatness of
$r(t)$.

 \proclaim{Lemma 4.4} Consider the $r(t,n)$ in Theorem 4.1. Then
there exist  $r_\pm \in \ell ^2$ such that $\|  r_\pm \| _{\ell ^2}
\le C \epsilon $ for fixed $C=C(\omega _0)$ and
$$\lim _{t\to \pm \infty }\| r(t)- e^{-it \Delta}r_{\pm}\| _{\ell ^2}=0.$$
\endproclaim
{\it Proof.}  We first write (4.2) as
$$\aligned & e^{ iHt}P_c(H)w(t)= P_c(H)w(0)+
\sum _{j=2}^4 e^{ iHt}w_j(t),
 \\&
e^{ iHt}w_j(t)=-i \int _0^t e^{ iHs } e^{ i\Theta (s)}P_c(H) (s)
g_j(s9ds.\endaligned $$ Then we observe that for $t_2>t_1$ by
(4.3-7)
$$\aligned \left \|   e^{ iHt_2}w_j(t_2)-e^{ iHt_1}w_j(t_1) \right \|
_{\ell ^2}\le \| g_j \| _{X_{(r,p)} (t_1,t_2)+\ell ^{2,2}
(t_1,t_2)}\to 0 \text{ for $t_1\to \infty$}.
\endaligned $$
This implies that the following limits exist  $$\lim _{t\to
+\infty}e^{ iHt}P_c(H)w(t)= \lim _{t\to +\infty}e^{ iHt} w(t)=w_+,
\tag 1$$ with the first equality due  to $w(t) \in C^1(\Bbb R, \ell
^2) \cap \ell ^{\frac 32 r}( \Bbb {Z},L^{\infty}_{t}([n,n+1],
    \ell ^{ p}(\Bbb {Z})))$ for any admissible pair $(r,p)$,
    so that $\lim _{|t|\to \infty}P_d(H) w(t)= 0$.
Recall that $u(t)= e^{i\Theta (t)}\phi _{\omega (t)} +w(t).$ Then
Theorem 4.1 and (1) imply
$$\aligned &
\lim_{t\to\infty}\|u(t )-e^{i\Theta(t)}\phi_{\omega_{+}} -e^{ -itH
}w_+\|_{\ell ^2}=0 .
 \endaligned $$
By Pearson's Theorem, see Theorem XI.7\cite{RS}, the following two
limits exist in $\ell ^2$,  for $w\in \ell^2_c(H )$ and $u\in \ell^2
$:

$$W u= \lim _{t\to +\infty }   e^{  it H
 }e^{  it\Delta }u \, , \quad    Z w= \lim _{t\to +\infty }  e^{ -it\Delta } e^{ -it H
 }w. $$
This follows from the fact that $H+\Delta =q$ with the operator
$u(n)\to q(n) u(n)$ in the trace class because of $q\in \ell ^1$.
For $u_+=Zw_+$ we have  the following, which yields Lemma 4.4:
$$  \lim_{t\to\infty} e^{ -itH } w_+  = \lim_{t\to\infty}
e^{  it\Delta  }u_+ \text{  in   $\ell ^2$}.$$

\bigskip
 \head \S 5 Dispersive theory for $H$ \endhead

We recall basic facts concerning the resolvent of the difference
Laplace operator. First,  for $g(\theta ) \in L^2(-\pi , \pi )$ and
for $  {u}(n )=\widehat{g}(n)$ we have

$$   - ({\Delta u})(n)=2 \left [ (1-\cos \theta
)g(\theta)\right ] ^{\wedge} (n).
$$
The resolvent $R_{-\Delta}(z)$  for $ z\in \Bbb C \backslash [0,4]$
has kernel
$$
R(m, n, z)=\frac{-i}{2\sin \theta }e^{-i \theta   |n-m|}, \ \ m,n
\in \Bbb Z,
$$
with $\theta$ the  unique solution     to $    2(1-\cos \theta ) =z
 $ in
$$  D=\{ \theta : -\pi \le \Re \theta \le \pi , \, \Im \theta <0 \}. $$
For all the above see \cite{KKK}. Then $(-\Delta -z)\psi=f, \ f\in
\ell^2(\Bbb Z)$ has solution
$$\psi_n =(R_{- \Delta}(z)f)_n
= \frac{-i}{2\sin \theta}\sum _{m\in \Bbb Z}e^{-i \theta
|n-m|}f_m.$$ We consider now discrete Jost functions. For $z$ and
$\theta$ as above we look for  functions  ${f}_{ \pm } (n, \theta )$ with
$$H{f}_{  \pm } (n, \theta )=z{f}_{  \pm } (n, \theta ) \text{ with }
\lim _{n\to \pm \infty}\left [ {f}_{  \pm } (n, \theta ) - e^{\mp in
\theta }\right  ]=0.\tag 5.1$$  Therefore the Green representation
of the solutions is
$$\aligned   {f}_{  \pm } ( n , \theta )=&
 e^{\mp in\theta }  - \sum _{m=n}^{\pm \infty
} \frac{\sin  (\theta (n-m))}{\sin \theta } q(m)  {f}_{ \pm } ( m,
\theta )  .\endaligned \tag 5.2
$$
Let  $
 {m}_{  \pm} $ be defined by  $ {f}_{  \pm } (n, \theta )=
 e^{\mp i n\theta
  } {m}_{  \pm } (n,\theta ).$
  Then we have
$$\aligned   {m}_{  \pm } (n,\theta )=&
 1  + \sum _{\nu =n}^{\pm \infty
} \frac{1-  e^{ 2i(n-\nu )\theta } }{2i\sin \theta   } q(\nu ) {m}_{
\pm } (\nu ,\theta )
 .\endaligned \tag 5.3
$$
Following standard arguments, see Lemma 1  \cite{DT}, we have:

 \proclaim{Lemma  5.1}
 For $\theta \in \overline{\Bbb C_-},$  (5.2) has for
any choice of sign a unique solution satisfying the estimates listed
below. These solutions solve $ H  u=zu$   with the asymptotic
property $ {f}_{ \pm } (n, \theta )\approx e^{\mp in \theta }
+o(e^{\mp in\theta})$ for $n\to \pm \infty .$ For $q\in \ell
^{1,\sigma}, \ \sigma\in [1,2],$
   there
 is a
  $C=C(q)$    such that
$\forall \, n\in \Bbb N$   we have,

$$  \align &
| {m}_{  \pm } (n,\theta)-  1 |\le C   \langle n^{\pm} \rangle
^{-\sigma} | \sin \theta |^{-1} e^{\frac{C}{| \sin \theta |}}, \
\theta \not \in \pi \Bbb Z, \tag 1
\\& | {m}_{  \pm }
(n,\theta)-  1 |\le C \langle n^{\pm} \rangle ^{-(\sigma-1)} \langle
\sin \theta  \rangle ^{-1} (1+n^{\mp} ) . \tag 2
\endalign
$$
${m}_{  \pm } (n,\theta)$ are for any $n$ analytic for $\theta \in
\Bbb C_-$, they satisfy ${m}_{  \pm } (n,\theta)={m}_{  \pm }
(n,\theta +2\pi )$, and extend into continuous functions in
$\overline{\Bbb C_-}$.

If $\sigma =2$,  there
 is a
  $C=C(q)$  such that

$$  \align &  \big |  \dot  {m}_{ + } (n,\theta)\big |\le C
  \langle
n^{-} \rangle ^2, \tag 3\\ &  \big |  \dot {m}_{ - } (n,\theta)\big
|\le C  \langle n^{+} \rangle ^2,   \tag 4
\endalign
$$
where $\dot {m}_{  \pm } (n,\theta)$ are for any $n$ analytic for
$\theta \in\Bbb C_-$, and extend into continuous functions in
$\overline{\Bbb C_-}$.
\endproclaim
  {\it Proof.} It is not restrictive to consider    $m _+(n,\theta
  ).$  We set $m(n,\theta )=m _+(n,\theta
  )$  in the rest of this lemma. For $D  (n , \theta ):
  =\frac{1- e^{ 2in\theta }   }{2i\sin \theta   }  $
we have
  $$ m(n,\theta )=1+\sum _{\nu  =n}^\infty D  (n-\nu  , \theta )q
  (\nu )  m(\nu ,\theta ).
\tag 5.4$$ We search for a solution

  $$ m(n,\theta )=1+\sum _{\ell =1}^\infty g_\ell (n,\theta ),\tag 5
$$
defined recursively by

$$g_\ell (n,\theta )=\sum _{n\le n_1 \le ...\le n_\ell}
D(n-n_1 , \theta ) \cdots D(n_{\ell -1} -n_\ell  , \theta
)q(n_1)\cdots q(n_\ell ).
$$
Notice that the $g_\ell (n,\theta )$ are $2\pi$ periodic in
$\theta$. Since this is true also for the estimates, we can assume
below that $\theta \in \overline{D}$. By  $|D(\mu , \theta )|\le
1/|\sin \theta |$   for $\mu \le 0$ we get
$$ \aligned & |g_\ell (n,\theta )|\le
\frac{1}{|\sin \theta |^\ell }\sum _{n\le n_1 \le ...\le n_\ell}
|q(n_1)\cdots q(n_\ell )|=  \frac{(\sum _{m=n}^\infty
|q(m)|)^\ell}{\ell ! |\sin \theta |^\ell}.\endaligned
$$
Therefore we get the following which yields (1)
$$ | m(n,\theta )-1|\le  \frac{1
}{|\sin \theta |  }    \eta (n) e^{ \frac{\eta (n) }{|\sin \theta |
} }.
$$ We consider  now
inequality (2). It is enough to assume $\theta \in \overline{D}$ is
close either to 0 or to $\pm \pi$, since otherwise estimate (1) is
stronger than (2). Then $ |D(n, \theta )|\le C_0|n|$ for a fixed
$C_0$. Then we get
$$ \aligned & |g_\ell (n,\theta )|\le C_0^\ell \sum _{n\le n_1 \le ...\le
n_\ell}    ( n_1-n) \cdots (n_\ell -n_{\ell -1}) | q(n_1)\cdots
q(n_\ell )|= \\& =  C_0^\ell \frac{(\sum _{m=n}^\infty (m-n)
|q(m)|)^\ell}{\ell !},\endaligned
$$
which yields

$$ | m(n,\theta )-1|\le C_0\gamma (n) e^{C_0\gamma (n)} .
$$
Notice that the bound increases exponentially if $n\rightarrow
-\infty.$ The following chains of inequalities are fulfilled

$$\aligned & | m(n,\theta ) | \le 1 +\sum _{ \nu =n}^\infty (\nu
-n) |q(\nu )| | m(\nu ,\theta ) | \le \\& \le 1 +\sum _{\nu
=n}^\infty \nu |q(\nu )| | m(\nu ,\theta ) | -n\sum _{\nu =n}^\infty
|q(\nu )| | m(\nu ,\theta ) |
\\& \le  1 +\sum _{\nu =0}^\infty \nu  |q(\nu )| |
m(\nu ,\theta ) | -n\sum _{\nu =n}^\infty   |q(\nu )| | m(\nu
,\theta ) |.
\endaligned$$
Furthermore we have also
$$1+\sum _{\nu =0}^\infty \nu  |q(\nu )| |
m(\nu ,\theta ) |\le 1+ \gamma (0)e^{\gamma (0)}
  \sum _{\nu
=0}^\infty \nu  |q(\nu )| =:K<\infty .$$ If we set

$$M(n,\theta ):=\frac{m(n,\theta )}{K(1+|n|)},$$
we   get the inequality
$$|M(n,\theta )|\le 1+
\sum _{\nu =n}^\infty  (1+|\nu |) |q(\nu )| | M(\nu ,\theta ) |,$$
 which can be solved by an iteration argument as in the previous case
  obtaining

$$|M(n,\theta )|\le  e^{
\sum _{\nu =n}^\infty  (1+|\nu |) |q(\nu )| | M(\nu ,\theta ) |}\le
K_1<\infty,
$$
and this gives
$$|m(n,\theta )|\le K_2 (1+|n|).  $$
Thus we get (2) by

$$\aligned & | m(n,\theta )-1 | \le  \sum _{\nu =0}^\infty \nu  |q(\nu )| |
m(\nu ,\theta ) | -n\sum _{\nu =n}^\infty   |q(\nu )| | m(\nu
,\theta ) | \le
\\& \le \gamma (0)e^{\gamma (0)}
\sum _{\nu =0}^\infty \nu |q(\nu)| -K_2 n\sum _{\nu =n}^\infty
 (1+|\nu |) |q(\nu )| | \le \\& \le K_3(1+\max (-n, 0))\sum _{\nu =n}^\infty
 (1+|\nu |) |q(\nu )| .
\endaligned$$
Notice that the above arguments yields also the uniqueness of
$m(n,\theta ).$  The $m( n ,\theta )$ are defined for all $\theta $
with $\Im \theta \le 0$. They satisfy $m( n ,\theta )=m( n ,\theta
+2\pi )$ and are analytic for $\Im \theta <0$. The last two
properties follow from the fact that the    $g(n ,\theta )$  satisfy
these properties and  that  series (5)  converges uniformly   for
$\Im \theta \le 0  $.

We now prove (3)   (the proof of (4) is similar). Better estimates
than (3) can be obtained using (1), the analyticity of $m( n ,\theta
)$ for $\Im \theta \le 0  $ and the Cauchy integral formula. So  it
is enough to assume that $\theta \in \overline{D}$ is close to the
interval $[-\pi , \pi ]$. Differentiating (5.4) we get
$$\dot m(n,\theta )=\sum _{\nu =n}^\infty D(n-\nu , \theta ) q(\nu )
\dot m(\nu ,\theta ) +\sum _{\nu =n}^\infty \dot D(n-\nu , \theta )
q(m) m(\nu ,\theta ).$$ We consider the representations
$$ \aligned & D(n-\nu , \theta )= \frac{\theta}{\sin \theta }\int ^{0}
_{ n-\nu }e^{ 2i\theta t}dt \text{ for  $\Re \theta \in [-\pi /2 ,
\pi /2 ] $},\\&    D(n-\nu , \theta )=\frac{\theta \mp \pi }{\sin
\theta }\int ^{0} _{n-\nu }e^{ 2i(\theta \mp \pi )t}dt \text{ for
$\Re \theta \in [ \pi /2 , \pi   ] $ (resp.   $\Re \theta \in [ -
\pi, \pi/2 ]$).}
\endaligned $$  Then

$$ \aligned & \dot D(n-\nu , \theta )=
\partial _\theta \left (
\frac{\theta}{\sin \theta }\right ) \int ^{0} _{n-\nu }e^{ 2i\theta
t}dt +2i\frac{\theta}{\sin \theta }\int ^{0} _{n-\nu }te^{ 2i\theta
t}dt,
\endaligned \tag5.5$$

$$ \align &
\dot D(n-\nu , \theta )=\partial _\theta \left (\frac{\theta \mp \pi
}{\sin \theta }\right ) \int ^{0} _{n-\nu }e^{ 2i(\theta \mp \pi
)t}dt +2i \frac{\theta \mp \pi }{\sin \theta } \int ^{0} _{n-\nu
}te^{ 2i(\theta \mp \pi )t}dt  . \tag 5.6
\endalign $$
By  $n-\nu \le 0$,  (5.5) implies for $\theta $ close to $[-\pi ,
\pi ]$
$$ | \theta \dot D(n-\nu , \theta )|\le C|n-\nu |,
\ \ \Re \theta \in [-\pi /2 , \pi /2 ]. $$ Similarly by   (5.6) we
obtain
$$ | (\theta - \pi ) \dot D(n-\nu , \theta )|\le C|n-\nu |, \ \ \Re \theta \in [  \pi/2 , \pi  ],$$
and
$$ | (\theta + \pi ) \dot D(n-\nu , \theta )|\le C|n-\nu |,\ \ \Re \theta \in [  -\pi , -\pi /2  ].$$
Furthermore we have the inequality
$$ |   \dot D(n-\nu, \theta )|
=\left |\partial_\theta\left (1-\frac{ e^{ 2i(n-\nu )\theta }
 }{2i\sin \theta   }\right )\right |\le C(n-\nu )^2,$$ and so
$$\aligned & \left |
\sum _{\nu =n}^\infty \dot D(n-\nu , \theta ) q(\nu ) m(\nu ,\theta
)\right | \le C \sum _{\nu =n}^\infty (n-\nu )^2 | q(\nu ) m(\nu
,\theta )| .
\endaligned$$
Suppose $n<0$.  Then from the fact that $q\in \ell ^{1,2}$

$$\aligned &  \sum _{\nu =n}^\infty \nu ^2 | q(\nu ) m(\nu ,\theta )| =
\sum _{\nu =n}^0 \nu ^2 | q(\nu ) m(\nu ,\theta )| +\sum _{\nu
=0}^\infty \nu ^2 | q(\nu ) m(\nu ,\theta )| \\& \le \sum _{\nu
=n}^0 n^2 | q(\nu ) m(\nu ,\theta )| +\sum _{\nu =0}^\infty \nu ^2 |
q(\nu ) m(\nu ,\theta )| \le K(1+n^2)
\endaligned$$
where we used $|m(\nu ,\theta )|\le C \langle \nu ^-\rangle $. For
$n\ge 0 $
$$\aligned &  \sum _{\nu =n}^\infty \nu ^2 | q(\nu ) m(\nu ,\theta )| \le K
\sum _{\nu =0}^\infty \nu ^2 | q(\nu ) |    .
\endaligned$$
Hence, for all $n$ we obtain

$$\aligned &  \sum _{\nu =n}^\infty \nu ^2 | q(\nu ) m(\nu ,\theta )| \le K
 \langle n^-\rangle ^2   .
\endaligned$$
For $n\ge 0$
$$\aligned &
\sum _{\nu =n}^\infty (n-\nu )^2 | q(\nu ) m(\nu ,\theta )| \le \sum
_{\nu =0}^\infty \nu ^2 | q(\nu ) m(\nu ,\theta )|\le K.
\endaligned$$
For $n< 0$ we obtain the chain of inequalities

$$\aligned &
\sum _{\nu =n}^\infty (n-\nu )^2 | q(\nu ) m(\nu ,\theta )| \le
2\sum _{\nu =n}^\infty \nu ^2 | q(\nu ) m(\nu ,\theta )| +  2
n^2\sum _{\nu =n}^\infty | q(\nu ) m(\nu ,\theta )| \\& \le K
(1+\nu^2) .
\endaligned$$
So
$$|\dot m(n,\theta )|\le \sum _{\nu =n}^\infty   (\nu -n) |q(\nu )
\dot m(\nu ,\theta ) |+K_2\langle n^-\rangle ^2,$$ and iterating
$$|\dot m(n,\theta )|\le K_2\langle n^-\rangle ^2 e^{\gamma (n)}.$$
We get, for any $n\in \Bbb N$,
$$|\dot m(n,\theta )|\le K_2\langle n^-\rangle ^2 +\sum _{\nu =n}^\infty   m |q(m)
\dot m(\nu ,\theta ) |-n\sum _{\nu =n}^\infty   |q(\nu ) \dot m(\nu
,\theta ) |.$$ The right hand side is smaller than
$$ \gamma (0)e^{\gamma (0)}  \sum _{\nu =0}^\infty m
|q(\nu )|-n\sum _{\nu =n}^\infty   |q(\nu ) \dot m(\nu ,\theta )
|,$$ which can be bounded by
$$
K_3\langle n^-\rangle ^2 (1+\sum _{\nu =n}^\infty
 \langle \nu \rangle ^2 |q(\nu )|).
$$
Following the same line of the proof of (1)-(2), by an iteration
argument we get the desired estimate and complete the proof of the
Lemma. Analyticity  of $\dot m(n,\theta )$ in the interior of $D$ and
continuity in $\overline{D}$ can be proved as the similar statement
for $  m(n,\theta )$.

\bigskip
\bigskip
For any fixed $n$ Lemma 5.1    implies    the Fourier expansion
$$m_{\pm }(n,\theta )=1+\sum _{\nu =1}^{\infty} B_{\pm }(n ,\nu )
e^{-i\nu \theta }.\tag 5.7$$ We consider the following:

\proclaim{Lemma 5.2} For $q\in \ell ^{1,1}$ we have
$$\sup _{n\ge 0} \| B_{+ }(n ,\nu ) \| _{\ell ^1_\nu
}<\infty  \text{ and } \sup _{n\le 0} \| B_{- }(n ,\nu ) \| _{\ell
^1_\nu }<\infty .
$$
\endproclaim
{\it Proof.} It is not restrictive to consider the $+$ case only. We
drop the $+$ subscript. By substituting $ e^{-in\theta }m(n,\theta
)=f_+(n,\theta )$ in (5.1) and using $z=2-2\cos (\theta )$ we obtain
$$e^{-i \theta } \left ( m(n+1,\theta
)-m(n,\theta ) \right )+e^{ i \theta } \left ( m(n-1,\theta
)-m(n,\theta ) \right )=q(n) m(n,\theta ).$$ Substituting the
Fourier expansion (5.7) we obtain
$$\aligned &  \sum _{\nu =1}^{\infty} \left ( B (n+1 ,\nu -1 )
 - B (n  ,\nu -1) \right )e^{-i\nu \theta }  +\\&
 \sum _{\nu =0}^{\infty}\left ( B (n-1 ,\nu +1) -B (n  ,\nu +1 )
 \right ) e^{-i\nu \theta }
=q(n)+q(n) \sum _{\nu =1}^{\infty} B (n ,\nu ) e^{-i\nu \theta }.
\endaligned
$$
   We have $B (n-1 , 1)- B (n , 1 )=q(n) $ and for  $\nu >0$
$$\aligned &
B (n-1 ,\nu +1)- B (n ,\nu +1 ) =q(n)B (n  ,\nu   )+ B (n   ,\nu
-1)- B (n +1  ,\nu -1 ) .
\endaligned
$$
By induction on $\nu $
$$\aligned & B (n  ,2\nu  )- B (n +1 ,2\nu   )=\sum _{j=1}^{\nu  }
q(n+j)B (n  +j,2(\nu -j)+1  )\\& B (n  ,2\nu -1 )- B (n +1 ,2\nu
-1)=q(n+\nu ) +\sum _{j=1}^{\nu  } q(n+j)B (n  +j,2(\nu -j)   ).
\endaligned
$$
So
$$\aligned & B (n  ,2\nu  ) =\sum _{l=n}^{\infty} \sum _{j=1}^{\nu  }
q(l+j)B (l  +j,2(\nu -j)+1  )\\& B (n  ,2\nu -1 ) = \sum
_{l=n+\nu}^{ \infty}q(l) +\sum _{l=n}^{\infty}\sum _{j=1}^{\nu  }
q(l+j)B (l +j,2(\nu -j) ).
\endaligned
$$
which after a change of variables, and setting $B_+(n,0)=0$, we
write as

$$\aligned & B (n  ,2\nu  ) =\sum _{l=0}^{\nu -1 }
\sum _{j=n+\nu -l}^{\infty} q(j) B (j,2l+1  )\\& B (n  ,2\nu -1 )
=\sum _{l=n+\nu}^{ \infty}q(l)  +\sum _{l=0}^{\nu -1 }\sum _{j=n+\nu
-l}^{\infty} q(j) B (j,2l   ).
\endaligned \tag 1
$$
By Lemma 5.1 (1) we know that $B\in \ell ^\infty ( \Bbb Z _{\ge
0}^2).$  We show now that (1) admits just one solution in this
space, which satisfies the bounds in the statement. We consider
$\Cal  B (n  ,\nu  )=\sum _{m =0}^{\infty}  K_m (n  ,\nu  )$ with
$$\aligned &
K_0 (n  ,2\nu -1 )=\sum _{l=n+\nu}^{ \infty}q(l) \text{ and }  K_0
(n  ,2\nu  )=0
\\& K_m (n  ,2\nu  ) =\sum _{l=0}^{\nu -1 }
\sum _{j=n+\nu -l}^{\infty} q(j) K_{m-1} (j,2l+1  )\\& K_m  (n ,2\nu
-1 ) = \sum _{l=0}^{\nu -1 }\sum _{j=n+\nu -l}^{\infty} q(j) K_{m-1}
(j,2l )
\endaligned \tag 2
$$
By the same inductive argument of p.139 \cite{DT} one can prove
$$| K_m (n  ,\nu  )|\le  \frac{\gamma  ^m (n)}{m!}\eta (n+[\nu  /2] ).
\tag 3$$ Indeed (3) is true for $m=0$. Assume (3) true for $m .$
Then we can write
$$
\aligned &|K_{m+1} (n  ,\nu  )|\leq \sum_{l=0}^{ [\nu /2] }\, \,
\sum _{j=n+[\nu /2] -l}^{\infty }|q( j)| \frac{\gamma  ^m
(j)}{m!}\eta (j+[l/2] )\leq
\\& \leq \eta (n+[\nu /2])
 \sum_{l=0}^{ [\nu /2] }\, \,
\sum _{j=n+[\nu /2] -l}^{\infty }|q( j)| \frac{\gamma  ^m
(j)}{m!}\leq \\ &  \le \eta (n+[\nu /2]) \left (  \sum_{l=0}^{ [\nu
/2] } \sum _{j=n+[\nu /2] -l}^{n+[\nu /2] }|q( j)| \frac{\gamma ^m
(j)}{m!} + \sum _{j=n+[\nu /2]}^{\infty }|q( j)| \frac{\gamma ^m
(j)}{m!} [\nu /2]  \right)
\\& \leq \eta (n+[\nu /2]) \left
(  \sum _{j=n }^{n+[\nu /2] }|q( j)| \frac{\gamma  ^m (j)}{m!}(j-n )
+ \sum _{j=n+\nu}^{\infty }| \, q( j)| \frac{\gamma ^m (j)}{m!} j
\right)\leq
\\ & \leq \eta (n+[\nu /2]) \sum _{j=n+1}^{\infty}|q( j)| \frac{\gamma
^m (j)}{m!} j = \frac{\gamma  ^{m+1} (n)}{(m+1)!}\, \eta (n+[\nu
/2]).
\endaligned
$$
Thus we have  $|\Cal B(n, \nu)|\leq e^{\gamma(n)}\eta (n+[\nu /2])
$. Therefore
$$
\| \Cal B(n ,\nu ) \| _{\ell ^\infty_\nu} \le e^{\gamma (n)}\eta (n)
\text{ and } \| \Cal B(n ,\nu ) \| _{\ell ^1_\nu }\le e^{\gamma
(n)}\sum_{\nu =0}^{\infty} \eta(n+\nu)\le e^{\gamma (n)}\gamma (n).
$$
Hence  $\Cal  B (n  ,\nu  )$ satisfies bounds as in the statement.
$U(n, \nu ):=B (n  ,\nu  )-\Cal B (n  ,\nu  )$ satisfies the
equation

$$\aligned &
\\& U(n  ,2\nu  ) =\sum _{l=0}^{\nu -1 }
\sum _{j=n+\nu -l}^{\infty} q(j) U (j,2l+1  )\\& U (n ,2\nu -1 ) =
\sum _{l=0}^{\nu -1 }\sum _{j=n+\nu -l}^{\infty} q(j) U (j,2l ).
\endaligned
$$
Iterating the above procedure we conclude $U(n ,\nu )= 0  $.

\bigskip

Given two functions $u(n)$ and $v(n)$ we denote by $[u,v](n)=
u(n+1)v(n)-u(n) v(n+1)$ the Wronskian of the pair $(u,v)$. If $u$
and $v$ are
  solutions of $Hw=zw$    then $[u,v]$ is constant.
Since the equations $Hu=\lambda u$ cannot have all solutions bounded
near $+\infty$ we have the following: \proclaim{Lemma 5.3} Let $q\in
\ell ^{1,1}.$ Then we have:

{\item {(1)}} If for a $\theta _0\in \{ 0, \pi ,-\pi \}$ we have
$f_{+}(n,\theta _0)\in \ell ^\infty $, then for $W(\theta ):=[f_{+}(
\theta  ), f_{-}( \theta  )]$ we have $W(\theta _{0})=0$. We will
call generic an $H$ such that $ W(\theta _{0})\neq 0 $ for all
$\theta _0\in \{ 0, \pi ,-\pi \}$.

{\item {(2)}} No element   $\lambda \in [0,4]$ can be an eigenvalue
of $H$ in $\ell ^2$.

{\item {(3)}} If   $\lambda \in \Bbb R$  is an eigenvalue of $H$,
then $\dim (H-\lambda )=1$.

\endproclaim
{\it Proof.} If $q\in \ell ^{1,2}$ by Lemma 5.1 we have $$ \lim
_{n\to \pm \infty}\|  1-m_\pm (n,\theta )\| _{L^\infty _\theta }\to
0.\tag 3$$ By the fact that the $ m_\pm (n,\theta )$ depend
continuously on $q\in \ell ^{1,1}$, which can be proved using the
arguments in Lemma 5.1, (3) is valid also for $q\in \ell ^{1,1}.$
This and the continuity in $\theta$ implies that for both signs
$m_\pm (n,\theta )\not \equiv 0$ as a function of $n$ for any fixed
$\theta$. If any of the three claims (1-3) is wrong, then for some
$\lambda \in \Bbb R$ all the solutions of $$(H-\lambda ) u=0\tag 4$$
are in $\ell ^\infty$. Let now $u_1 (n)$ and  $u_2 (n)$ be a
fundamental set of such solutions. Consider now the equation
$(-\Delta   -\lambda ) U=0$ which we rewrite as $(H-\lambda ) U=qU.$
Then solutions $U\in \ell ^ \infty ([N,\infty ))$ can be written for
$q\in \ell ^1$  as
$$U(n)=u (n)-\sum _{j=n}^\infty  \frac{u_1(n)u_2(j)-u_2(n)u_1(j)}{
[u_1,u_2]}q(j) U(j)\tag 5$$ with $u (n)$ a solution of (4). But for
$q\in \ell ^1$ and $N$ large, (5) establishes an isomorphism inside
$\ell ^ \infty ([N,\infty ))$ between solutions of (4), which form a
2 dimensional space, and of $(-\Delta   -\lambda ) U=0$, which form
a 1 dimensional space. Obviously this is absurd. Therefore it is not
possible for all solutions of (4) to be in $\ell ^\infty$.

Next we introduce transmission  and reflection coefficients.
\proclaim{Lemma 5.4} Let $q\in \ell ^{1,1}$. For $\theta \in [-\pi ,
\pi ]$ we have $\overline{f_{\pm}(n,\theta )} =f_{\pm}(n,-\theta )$
and for $\theta \neq 0,\pm \pi $ we have
$$ {f}_{  \mp } (n,\theta )=
\frac{1}{ T  (\theta )}\overline{ {f}_{  \pm } (n,\theta )}
+\frac{R_{ \pm} (\theta)}{T  (\theta  )}  {f}_{  \pm } (n,\theta )
\tag 1 $$ where $T  (\theta )$ and $R_\pm (\theta  )$ are defined by
(1) and satisfy: $$  \align   & [\overline{ {f}_{ \pm } (\theta )},
{f}_{ \pm } (\theta )] = \pm 2i\sin \theta , \tag 2\\& T  (\theta )=
\frac {\pm 2 i\sin \theta } {[ f_{ \mp } (\theta ),
 f_{ \pm } (\theta )]  } \, , \quad R_{  \pm } (\theta ) = -\frac {[
 {f_{ \mp } ( \theta  )},  \overline{f_{  \pm }} (\theta)] }
{[ f_{ \mp } ( \theta ),  f_{  \pm } (\theta )] } ,\tag 3
\\&
  \overline{T  (\theta  )}=T  (-\theta  ) \, , \,
\overline{R_{  \pm } (\theta  )}=R_{  \pm } (-\theta  ) , \tag 4
\\&|T  (\theta  )|^2+ |R_{ \pm } (\theta ) |^2=1\, ,
\quad T (\theta )  \overline{R_{ \pm } (\theta  )}+R_{  \mp} (\theta
) \overline{ T  (\theta  )}=0.\tag 5\endalign
$$
\endproclaim

{\it Proof.}   $\overline{f_{\pm}(n,\theta )} =f_{\pm}(n,-\theta )$
follows  by the fact that $q(n)$ has real entries and by uniqueness
in Lemma 5.1. The pair in the right hand side of (1) is linearly
independent, so (1) follows from the properties of second order
homogeneous linear difference equations. (2-4) follow applying
Wronskians to (1).   Iterating (1)  twice we get (5). Indeed, for
example

$$\aligned & f_-= \frac{1}{T}\overline{f_+}+\frac{R_+}{T}f_+=
\frac{1}{T}\left( \frac{1}{\overline{T}}
{f_-}+\overline{\frac{R_-}{T}f_-} \right ) +\frac{R_+}{T}\left(
\frac{1}{ {T}} \overline{{f_-}}+ \frac{R_-}{T}f_-  \right ) \\& =
\left( \frac{1}{ |T|^2}  + \frac{R_+R_-}{T^2} \right ) f_- + \left(
\frac{\overline{R_-}}{ |T|^2}  + \frac{R_+ }{T } \right )
\overline{f_-}.
\endaligned $$
This yields $R_-\overline{T}+\overline{R_+}T=0$ and $ \frac{1}{
|T|^2}  + \frac{R_+R_-}{T^2}=1$.  Substituting
$R_-=-\frac{\overline{R_+}T}{\overline{T}}$ we get
$|T|^2+|R_+|^2=1$. Similarly one gets $|T|^2+|R_-|^2=1$.

\bigskip

 \proclaim{Lemma 5.5} Let $W(\theta ):=[
 f_+(\theta ), f_-(\theta )]$ and $W_1(\theta ):=[
 f_+(\theta ), \overline{f}_-(\theta )]$.

 {\item {(1)}}   For $\theta \in [-\pi ,\pi ]\backslash
 \{ 0, \pm \pi \}$ we have $W(\theta )\neq 0$.
 We have $|W(\theta )|\ge 2|\sin \theta |$ for all
 $\theta \in [-\pi , \pi ]$
 and  in the
generic case $|W(\theta )|>0$.

{\item {(2)}} For  $j=0,1$ and  $q\in \ell ^{1,1 +j}$ then $W(\theta
)$ and $W_1(\theta )$  are in  $  C^j[-\pi , \pi ]$.

{\item {(3)}}  If $q\in \ell ^{1,2}$ and $ W(\theta _{0})=0$ for a $
\theta _{0}\in \{ 0, \pm \pi \}$, then $\dot W(\theta _{0})\neq 0$.
In particular if $q\in \ell ^{1,2}$, then $T  (\theta  )=  - { 2
i\sin \theta }/ {W(\theta )  }$  can be extended continuously in
$[-\pi , \pi ]$ with $T(-\pi )=T(\pi )$.

\endproclaim
{\it Proof}. (1) follows immediately from   $T  (\theta  )=  - { 2
i\sin \theta }/ {W(\theta )  }$ and $|T  (\theta  )|\le 1$ and the
definition of $H$ generic  in  Lemma 5.3. (2) follows from Lemma
5.1. (3) follows from Lemma 5.1 and (1).

\bigskip
We need of the following: \proclaim{Lemma 5.6} Let   $0>\Im \theta $
and $q\in \ell ^{1,1}$ with $z=2(1-\cos \theta )$ not an eigenvalue
of $H$. Then the resolvent $R_{H } (z)$ has kernel $R_{H } (n,m ,z)$
such that $R_{H } (n,m ,z)=K  (n,m  )$ where we define

$$\aligned &  K  (n,m  ) =-\frac{f_{ - } (n,\theta )
   f_{ + } (m,\theta )  }
{[ f_{ + } (\theta ),   f_{ - } (\theta ) ] } \text{  for $n<m$} \\
& K  (n,m  )=-\frac{f_{ + } (n,\theta )f_{ - } (m,\theta )
     }
{[ f_{ + } (\theta ),   f_{ - } (\theta ) ] } \text{  for $n\ge m$}.
\endaligned $$
\endproclaim
{\it Proof.} First of all, if $z\not \in [0,4]$ is not an eigenvalue
of $H$ we have $[ f_{ + } (\theta ),   f_{ - } (\theta ) ] \neq 0$
and $K(n,m)$ is well defined. Indeed if $[ f_{ + } (\theta ),   f_{
- } (\theta ) ]$ then $f_\pm (\cdot , \theta )$ are  proportional.
Then they   belong to $\ell ^2$ and in particular are eigenvectors
with eigenvalue $z$.
 Furthermore, by (1) Lemma 5.1 we have  $|K(n,m)|\le C(\theta )
e^{|n-m|\Im (\theta )}$ for some constant $ C(\theta )$. Then
$K(n,m)$ is the kernel of an operator $K\in B(\ell ^2, \ell ^2)$.
For any fixed $m$ by definition of $H$ we have
$$HA (\cdot ,m  ) (n)= 2K (n,m  )-K(n+1,m  )-K (n-1,m  )+q(n)KA (n,m  )
 .$$ By elementary verification for $n>m$ from the above
identity we get $$(H-z)K (\cdot ,m  ) (n)=-((H-z)f_+) (n)  f_-(m)
/W(\theta) =0,$$ while for $n<m$ we obtain similarly $$(H-z)K (\cdot
,m ) (n)=-((H-z)f_-) (n)  f_+(m) /W(\theta) =0.$$  Finally, in the
$n=m$ case

$$\aligned & (H-z)K (\cdot ,m  ) (m)=-((H-z)f_+) (m)
 \frac{f_-(m)}{   W(\theta)}  \\&+
\frac{f_+ (m) f_- (m-1)-f_+ (m-1) f_- (m )}{{
W(\theta)}}=1.\endaligned $$ This implies that
$(H-z)K=1=(H-z)R_H(z)$ and so   $K (n,m  ) =R_{H } (n,m ,z) $. This
concludes Lemma 5.6.

\bigskip

Next we have, see also Theorems 1 and 2 \cite{PS}:

\proclaim{Lemma 5.7} Let $q\in \ell ^{1,1}(\Bbb Z)$.

{\item {(a)}} Assume   $H$  is generic in the sense of Lemma 5.3.
Then for $\sigma
>1$ we have that for $\lambda \in \sigma _c(H)$ the  following
 limit exists   in $C^0([0,4],B(\ell ^{2,\sigma } ,
\ell ^{2,-\sigma }    ))$
$$\lim _{\epsilon \to 0 ^+} R_{H }(\lambda \pm i\epsilon )
=R^{\pm}_{H }(\lambda  ) . \tag 1$$ Furthermore, for $\lambda $ is
some fixed small neighborhood  in $\Bbb R$ of $\sigma _c(H)$, and
for $\epsilon >0$, the operators $R_{H }(\lambda \pm i\epsilon )$
are Hilbert Schmidt  (H-S) with H-S norm uniformly bounded.

{\item {(b)}}If $H$  is  not generic, then (1) exists  pointiwise
for $\sigma
>1 $  in $B(\ell ^{2,\sigma } ,
\ell ^{2,-\sigma }    )$  for any $0<\lambda <  4  $.
\endproclaim
{\it Proof}. We start with (a). It is not restrictive to consider
the limit from above. For $n\ge m$, by  Lemma 5.6 we have for
$z_{\epsilon }=\lambda +i \epsilon $ and for the corresponding
$\theta _\epsilon$
$$\langle n \rangle ^{-\sigma} \langle m \rangle
 ^{-\sigma}R_{H } (n,m ,z_{\epsilon }) =-\langle n
 \rangle ^{-\sigma} \langle m \rangle ^{-\sigma}
 \frac{f_{ + } (n,\theta _\epsilon )f_{ - } (m,\theta _\epsilon )
     }
{[ f_{ + } (\theta _\epsilon ),   f_{ - } (\theta _\epsilon) ] }.$$
The fact that $H$ is generic implies
$$ [ f_{ + } (\theta _\epsilon ),   f_{ - } (\theta _\epsilon) ]
\ge C \ge 0$$ for a fixed $C$.
  Lemma 5.1 implies  for $n\ge m$ there is a fixed $C>0$ such that
  $|f_{ + } (n,\theta _\epsilon )f_{ - } (m,\theta _\epsilon )|
  \le C(  1+  \max ( -n, 0)+   \max ( m, 0) ).$
Then we conclude that for a fixed $C>0$ we have

$$ \aligned & \sum _{m\in \Bbb Z} \sum _{n =m}^{\infty}
\left | \langle n
 \rangle ^{-\sigma} \langle m \rangle ^{-\sigma}
 \frac{f_{ + } (n,\theta _\epsilon )f_{ - } (m,\theta _\epsilon )
     }
{[ f_{ + } (\theta _\epsilon ),   f_{ - } (\theta _\epsilon) ] }
\right |^2 \lesssim \\& \sum _{n\ge m} \langle n
 \rangle ^{-2\sigma} \langle m \rangle ^{-2\sigma}(  1+  \max ( -n, 0)+   \max ( m, 0) )^2
  = I_1+I_2+I_3, \endaligned $$
where $I_1$ involves the sum for $|n|\approx |m|$, $I_2$    for
$|n|\gg |m|$ and $I_3$    for $|n|\ll |m|$. We have for $j=1$ $$
I_j\lesssim \sum   \langle n
 \rangle ^{-2\sigma +1} \langle m \rangle ^{-2\sigma +1 }<\infty . \tag 2$$
We have $$ I_2  =\sum _{n\ge  m, |n|\gg |m|} \langle n
 \rangle ^{-2\sigma}
  \langle m \rangle ^{-2\sigma}
  (  1+  \max ( -n, 0)+   \max ( m, 0) )^2
.
$$
 But $n\ge m$ and $ |n|\gg |m|$  implies $n\ge 0$
 and  so we get (2) for $j=2$.
We have $$ I_3  =\sum _{n\ge  m, |m|\gg |n|} \langle n
 \rangle ^{-2\sigma}
  \langle m \rangle ^{-2\sigma}
  (  1+  \max ( -n, 0)+   \max ( m, 0) )^2
.
$$
But $n\ge m$ and $ |m|\gg |n|$  implies $m\le 0$.  So we get (2) for
$j=3$.

 By a similar argument
$$ \aligned & \sum _{n\in \Bbb Z} \sum _{m =n+1}^{\infty }
\left | \langle n \rangle ^{-\sigma} \langle m \rangle
 ^{-\sigma}R_{H } (n,m ,z_{\epsilon })
\right |^2 \le  C_\sigma
 <\infty .\endaligned $$
So $R_{H }(\lambda + i\epsilon )$ are H-S, with a uniform bound on
the corresponding H-S  norm. By the continuous dependence of $f_{
\pm } (n,\theta   )$ and of the Wronskian on $\theta$, it is
elementary to conclude that the limit (1) holds in $C^0([0,4],B(\ell
^{2,\sigma } , \ell ^{2,-\sigma }    ))$ in the generic case. By the
same arguments (1) holds for $\lambda \in (0,4)$ in the non  generic
case.

\bigskip Notice that Lemma 5.7 (a) implies the first two claims of
Lemma 3.2. The last claim of Lemma 3.2 follows from: \proclaim{Lemma
5.8} For any $u\in \ell ^2(\Bbb Z)$ we have
$$\aligned & P_c( H )u=\lim _{\epsilon \to 0^+}
 \frac 1{2\pi i}  \int _{ [0, 4  ]}
 \left [ R_{H  }(\lambda +i\epsilon
)- R_{H  }(\lambda -i\epsilon )
 \right ] ud\lambda   .
\endaligned \tag 1
$$ In particular, for $u\in \ell ^{1,\sigma }$ with $\sigma >1$, we
have $$\aligned & P_c( H )u=
 \frac 1{2\pi i}  \int _{ [0, 4  ]}
 \left [ R_{H  }^+(\lambda
)- R_{H  }^-(\lambda   )
 \right ] ud\lambda   .
\endaligned \tag 2
$$
\endproclaim
{\it Proof.} (1) a consequence of the
 spectral theorem, see p.81 \cite{T}, while (2) holds because
the right hand side of (1) converges in $\ell ^{2,-\sigma}$ to the
right hand side of (2).

\proclaim{Lemma 5.9} Let $q\in \ell ^{1,1}$ if $H$ is generic and
$q\in \ell ^{1,2}$ if $H$ is not generic. For $u\in \Cal S(\Bbb Z)$
the following are well defined: $$\aligned & \int _{ 0}^4
   R_{H  }^{\pm}(\lambda
) ud\lambda  = \int _{ 0 }^\pi
  \sum _{\nu =-\infty }^{\infty } K_{\pm } (n, \nu ,\theta )
 u (\nu )  \sin \theta \,  d\theta   \\&
 = \sum _{\nu =-\infty }^{\infty }\int _{ 0 }^\pi
   K_{\pm } (n, \nu ,\theta ) \sin \theta \,  d\theta \, u (\nu ),
   \endaligned \tag 1$$
with  $$ \aligned & K_\pm (n,\nu ,\theta )= -\frac{ {f}_+(n,\pm
\theta )\, {f}_- (\nu ,\pm \theta ) }{W(\pm \theta )} \quad
\text{for $n\ge \nu  $ ,}
\\ & K_\pm (n,\nu ,\theta )=-\frac{ {f}_-(n,\pm \theta )
 {f}_+ (\nu ,\pm \theta )  }{W(\pm \theta )} \quad \text{for $n<\nu $.}
\endaligned \tag 2
$$

\endproclaim
{\it Proof.} By  Lemmas 5.6-7 the kernel of $R^{\pm }_H(\lambda )$
is given by $ K_\pm $ for $2-2\cos \theta =\lambda$ with $\theta \in
[0,\pi ]$. The first equality  in (1) is then a a consequence of a
change of variables. The second equality in (1) is consequence of
 Fubini equalities. The required summability for
  $H$    generic  follows from $|W(\theta )|>C>0$
  and for $H$   non generic from the fact that
$\sin \theta /W(\pm \theta )$ is a continuous function by Lemma 5.5.

\bigskip
We now recall Theorem 1.3:

\proclaim{Theorem 5.10} We assume $q\in \ell ^{1,2}$ in the non
generic case and $q\in \ell ^{1,1}$ in the   generic case. Then we
have:

$$\| P_c(H)e^{itH} : \ell ^1(\Bbb Z) \to \ell ^\infty (\Bbb Z)\| \le
C \langle t \rangle ^{-1/3}  \text{ for a fixed $C>0$.}$$
\endproclaim
{\it Proof of generic case.} It is not restrictive here to
assume $n<\nu$. We have
$$\aligned  &  P_c(H)e^{itH}(n, \nu )= \frac{1}{2\pi i}
\int _{ 0 }^\pi e^{it  (2-2\cos \theta )}
   \left [ K_{+ } (n, \nu ,\theta ) -
 K_{- } (n, \nu ,\theta )\right ]  \sin \theta \,  d\theta \\& =
   \frac{1}{2\pi i}
\int _{ -\pi  }^\pi e^{it  (2-2\cos \theta )+i\theta (n -\nu )}
    K  (n, \nu ,\theta ) d \theta \\&
       K  (n, \nu , \theta ):=
 {m_-}  (n, \theta )   m_+  (\nu , \theta )  \frac {\sin ( \theta
) } {W(\theta )}.  \endaligned \tag 1
$$
For $\Cal M  (\Bbb Z)$ the space of complex measures in $\Bbb Z$ we
have for $(n,\nu )$ fixed and taking Fourier series in $\theta$
$$ \aligned & \left |   (1) \right |
\le \left \| \left [ e^{it  (2-2\cos \theta )}   \right ] ^{ \vee }
        \right \| _{\ell ^\infty   }
\left \| \left [ K (n, \nu , \theta ) \right ] ^{ \wedge }
        \right \| _{\Cal M  (\Bbb Z) }\\&  \le  C \langle t
     \rangle ^{-\frac{1}{3}} \left \|
     \left [ K (n, \nu , \theta ) \right ] ^{ \wedge  }
       \right \| _{\Cal M  (\Bbb Z) },
        \endaligned $$
with the second inequality due to stationary phase. We have
$W^{\wedge}\in \ell ^1$. This follows from  (5.7) and Lemma 5.2.
Since $W(\theta )\neq 0 $ for all $\theta \in [-\pi , \pi ]$,  then
$
     \left [ 1/{W(\theta )} \right ] ^{ \wedge } \in \ell ^1$
  by Wiener's Lemma, see 11.6 \cite{R}. By
convolutions, we have $
     \left [ \sin ( \theta
) /{W(\theta )} \right ] ^{ \wedge } \in \ell ^1$. If $n\le 0\le
\nu$ we exploit, for $\delta _{a,0}$ the Kronecker delta,
$$  \| \widehat{m}_{+ }(\nu  ,a ) -\delta _{a,0} \| _{\ell ^1_a
}+ \| \widehat{m}_{- }(n ,a )-\delta _{a,0} \| _{\ell ^1_a }\le C
$$
for a fixed $C$.  Then for a fixed $C  $ we have
$$ \left \|
     \left [ K(n, \nu , \theta ) \right ] ^{ \wedge }
      \right \| _{\Cal M  (\Bbb Z) }\le C<\infty, \tag 2$$
for all $n\le 0\le \nu$.  If $0<n\le \nu  $, we can substitute $m_{
-} (n, \theta )$ using

$$\aligned &  2i\sin (\theta )\, {m_{ -}} (n, \theta ) =- {W (\theta )}
 \overline {m_{ +}} (n, \theta )+   {W _1(\theta )}  e^{-2in\theta}  {m_{ +}} (n, \theta ).\endaligned
$$
and we can repeat the argument. The argument for $ n\le \nu <0 $ is
similar.
\bigskip {\it Proof of non generic case.} By Lemma 5.5, $
\sin (\theta )/W(\theta ) $ is continuous and periodic. Suppose now
that $W(0)=  0$  and $W(\pm \pi )\neq  0$ . We consider a smooth
partition of unity $1=\chi +\chi _1$ on $\Bbb T=\Bbb R/(2\pi \Bbb
Z)$ with $\chi
 =1$ near $0$ and $\chi
 =0$ near $\pi $. Then it is enough to consider (the case with
   $\chi _1 $ can be treated  as above)

$$ \aligned  &  \int _{  -\pi  } ^{  \pi  }
 e^{it  (2-2\cos \theta )
     +i(n-\nu )\theta }   K (n, \nu , \theta )    \frac{\chi (\theta)}{\widetilde{\chi }(\theta)}
     d \theta
      ,\endaligned \tag 3
$$
  $\widetilde{\chi }(\theta)$   another cutoff with
$\widetilde{\chi }=1$ on the support of $\chi$  and $\widetilde{\chi
}=0$ near $ \pm \pi$.   We set
$$\frac{1} {\sin (\theta )} = \frac{1} {2\tan (\theta /2)}+
\phi (\theta), \quad  \phi (\theta):=\frac{1} {\sin (\theta )} -
\frac{1} {2\tan (\theta /2)}.$$    Then $\widetilde{\chi}(\theta )
\phi (\theta)\in C^\infty (\Bbb T)$. Since $
   W   ^{\wedge} \in \ell ^1$,   it follows that also $\left [
\widetilde{\chi}  \phi   W    \right ] ^{\wedge} (n)\in \ell ^1$. If
$q\in \ell ^{1,2}$, we have   $
    W   ^{\wedge}  \in \ell  ^{1,1}.  $ Indeed for
$m_\pm (n,\theta )=1+\widetilde{m}_\pm (n,\theta )$ we have
     $$ \aligned & W(\theta ) =-2i \sin \theta + e^{-i\theta}
\left (\widetilde{m}_+ (n+1,\theta )+\widetilde{m}_- (n,\theta )+
\widetilde{m}_+ (n+1,\theta ) \widetilde{m}_- (n,\theta )\right )
-\\& - e^{ i\theta} \left (\widetilde{m}_+ (n ,\theta
)+\widetilde{m}_- (n+1,\theta )+ \widetilde{m}_+ (n ,\theta )
\widetilde{m}_- (n+1,\theta )\right ).\endaligned $$ Then
 $
    W   ^{\wedge}  \in \ell  ^{1,\sigma -1}   $ is a consequence of
    $ B_{\pm}(n, \nu ), B_{\pm}(n+1, \nu ) \in \ell  ^{1,\sigma -1},$
    for $\sigma =1,2.$ This last fact for $q\in \ell  ^{1,\sigma} $
     follows from $|B(n,\nu )|\le e^{\gamma (n)}
     \eta (n+[\nu /2]) \lesssim   e^{\gamma (0)}
     \eta (n+[\nu /2])
     $, proved in Lemma 5.2.
     Indeed
     $$\aligned & \sum _{\nu =1}^\infty
     \nu^{\sigma -1} |B(n,\nu )|
     \lesssim \sum _{\nu =1}^\infty \nu^{\sigma -1}
     \sum _{j=n+[\nu /2]  }^\infty
|q(j)| \le \sum _{j=n   }^\infty |q(j)| \sum _{\nu
=1}^{j-n+1}\nu^{\sigma -1}\lesssim   \| q\| _{\ell ^{1,\sigma }}.
\endaligned $$
For $B_{\pm}(n+1, \nu )$ the argument is the same. Having
established $
    W   ^{\wedge}  \in \ell  ^{1,1},$  we have, see for example p.3
    \cite{Ch},
$$\aligned \left [ \frac{\widetilde{\chi}(\theta ) W(\theta )}
{-i 2\tan (\theta /2)} \right ] ^{\wedge} (n) = \sum _{\nu >n}\left [
 \widetilde{\chi}  W  \right ]
^{\wedge} (\nu ) -\sum _{\nu <n}\left [
 \widetilde{\chi}  W  \right ]
^{\wedge} (\nu ).
 \endaligned $$
Since $\widetilde{\chi} (0) W(0)=0$ and since $\left [
 \widetilde{\chi}  W  \right ]
^{\wedge} \in \ell ^{1,1}$, it follows $\left [
\frac{\widetilde{\chi}(\theta ) W(\theta )} { \tan (\theta /2)}
\right ] ^{\wedge}\in \ell ^{1 }$. Hence we have proved $\left [
\frac{\widetilde{\chi}(\theta ) W(\theta )} { \sin (\theta  )}
\right ] ^{\wedge}\in \ell ^{1 }$. Then (3) can be bounded  with the
argument of the generic case. If $0<n\le \nu  $ we can show in a
similar way that $\left [ \frac{\widetilde{\chi}(\theta ) W_1(\theta
)} { \sin (\theta  )} \right ] ^{\wedge}\in \ell ^{1 }$. If also
$W(\pi )=0$ we can repeat a similar argument near $\pi$.

\bigskip

\head Appendix A: proof of Lemma 1.1 \endhead

Let us search for  a solution of (1.3) in the form $u=a\varphi
_0+h,$ with $h(n)\in \Bbb R$ for all $n$,
   $\langle h,\varphi _0\rangle _{\ell ^2}=0$ and $a\in \Bbb R$.  Then (1.3) becomes
$$\align  & (E_0-\omega )a+  a^7
    \|  \varphi  _0 \|^8_{\ell ^8} + \langle
N(h), \varphi _0 \rangle  =0\tag A.1
\\  h -R_H(-E_0)& P_c(H) \left
[ (E_0-\omega )h+    (a\varphi _0+h)^7 \right ] =0 ,\tag A.2
\endalign
$$
where  $  N(h) (n)= \sum _{j=1}^{7}  \left ( \matrix  7 \\ j
\endmatrix \right )(a \varphi _0 (n))^{7-j}
  \,  (h(n)) ^ j .$ We have $$R_H(-E_0)  P_c(H)\in
  B(\ell ^{p,\sigma }(\Bbb Z,\Bbb R), \ell ^{p,\sigma }(\Bbb Z,\Bbb R)) \text{
 for any $p\in [1,\infty ]$ and $\sigma \in \Bbb R$}.$$ The functions
in (A.1)--(A.2) are $C^\omega $ in the arguments $a,\omega \in \Bbb
R$ and $h\in \{ \varphi _0 \} ^{\perp}\cap \ell ^{p,\sigma }$.
Substituting  $h=a^7g$ and factoring out in (A2) we obtain for $g$
$$g -R_H(-E_0)  P_c(H) \left
[ (E_0-\omega )g+    ( \varphi _0+a^6g)^7 \right ] =0 .\tag A.3
$$
By the implicit function theorem applied to (A.3)  we have
$h=h(a,\omega )=a^7g(a^6,\omega )$ with $g(a^6,\omega )$ real
analytic in $(a^6,\omega )$  and   with values in $ \{ \varphi _0 \}
^{\perp}\cap \ell ^{p,\sigma }(\Bbb Z,\Bbb R)$.  Plugging in (A.1)
we obtain
$$\align  & (E_0-\omega ) +  a^6
    \|  \varphi  _0 \|^8_{\ell ^8} +\sum _{j=1}^{7}   \left ( \matrix  7 \\ j
\endmatrix \right ) a^{6(j+1)} \langle
 \varphi _0^{7-j}g^j(a^6,\omega ), \varphi _0 \rangle  =0.
\endalign
$$
By the implicit function theorem we obtain an analytic function $
\omega -E_0 \to a^6$ with
$$ a ^{6}=
    \|  \varphi  _0 \|^{-8}_{\ell ^8} ( \omega -E_0)
    (1+ O(\omega -E_0 )) .$$
 Then we obtain Lemma 1.1.

\bigskip

\head Appendix B: proof of Lemma 2.1 on global well posedness
\endhead
 The operator $iH$ is a bounded skew adjoint operator in $\ell ^2$
 and the nonlinearity $(F(u)) (n)=|u(t,n)|^ 6u(t,n)$ is Lipschitz
 continuous on bounded sets in $\ell ^2$. As a consequence
 we have what follows.

 {\item {(1)}} For any $u_0\in \ell ^2$ there exist $T_1(u_0)<0
 <T_2(u_0)$  and a solution $u(t)\in C^\infty ((T_1(u_0),T_2(u_0)), \ell ^2) $ of
  (1.2) with $u(0)=u_0$. If $T_j(u_0)\in \Bbb R$ for a $j$, then
$$\lim _{t\to T_j(u_0)}\| u(t)\| _{\ell ^2}=\infty .$$
For any solution $v(t) $ of the same Cauchy problem in $(\alpha
,\beta )\subset (T_1(u_0),T_2(u_0))$, then $v(t)=u(t)$ in $(\alpha
,\beta )$.  If $u_{0,\nu} \to u_0$ in $\ell ^2$ and for any bounded
interval $[a,b]\subset (T_1(u_0),T_2(u_0))$, then for $\nu $ large
the corresponding solutions $u_\nu (t)$ are in   $C^1([a,b], \ell
^2) $   and converge uniformly to $u(t)$ therein, see \cite{CH}
sections from 4.3.1 to 4.3.3. By $iu_t=-Hu-|u|^6u$ and by the fact
that the rhs is in $C^1(T_1(u_0),T_2(u_0))$   we conclude that $u\in
C^2(T_1(u_0),T_2(u_0))$, and so by induction $u\in
C^k(T_1(u_0),T_2(u_0))$ for all $k$.  The continuity with respect to
the initial data in $C^k_{loc}(T_1(u_0),T_2(u_0))$ is obtained
similarly from the continuity in
   $C^1_{loc}(T_1(u_0),T_2(u_0))$.

{\item {(2)}} For solutions $u(t)$ of (1.2) we have $\| u(t)\|
_{\ell ^2}=\| u(0)\| _{\ell ^2}.$ As a consequence, for any $u_0\in
\ell ^2$ we have $T_2(u_0)=+\infty$ and  $T_1(u_0)=-\infty$.

\enddocument